\documentclass[sn-mathphys,Numbered]{sn-jnl}% Math and Physical Sciences Reference Style
\usepackage{graphicx}%
\usepackage{multirow}%
\usepackage{amsmath,amssymb,amsfonts}%
\usepackage{amsthm}%
\usepackage{mathrsfs}%
\usepackage[title]{appendix}%
\usepackage{xcolor}%
\usepackage{textcomp}%
\usepackage{manyfoot}%
\usepackage{booktabs}%
\usepackage{algorithm}%
\usepackage{algorithmicx}%
\usepackage{algpseudocode}%
\usepackage{listings}%
\usepackage{chngcntr}
\theoremstyle{thmstyleone}%
\newtheorem{theorem}{Theorem}%  meant for continuous numbers
\newtheorem{proposition}{Proposition} % to get separate numbers for theorem and proposition etc.

\theoremstyle{definition}%
\newtheorem{example}{Example}%
\theoremstyle{definition}%
\newtheorem{definition}{Definition}%

\newcommand{\ntr}{\underline{\mathrm{tr}}}
\newcommand{\tr}{\mathrm{tr}}
\newcommand{\Tr}{\mathrm{Tr}}

\newcommand{\E}{\mathbb{E}}
\newcommand{\Ord}{\mathcal O}
\begin{document}

\title[Article Title]{
%Spectrum of subblocks of structured random matrices~: A free probability approach\\
Structured random matrices and cyclic cumulants~: A free probability approach}


\author[1]{\fnm{Denis} \sur{Bernard}}\email{denis.bernard@ens.fr}
\author*[1]{\fnm{Ludwig} \sur{Hruza}}\email{ludwig.hruza@ens.fr}

%\author[1]{{Denis} {Bernard}~\footnote{E-mail: {denis.bernard@ens.fr}.} }
%\author*[1]{{Ludwig} {Hruza}~\footnote{E-mail: {ludwig.hruza@ens.fr}, corresponding author.}}

\affil[1]{\orgdiv{Laboratoire de Physique} de \orgname{Ecole Normale Sup\'erieure, CNRS, ENS \& PSL University, Sorbonne Universit\'e, Universit\'e Paris Cit\'e}, \orgaddress{\street{24 rue Lhomond}, \city{Paris}, \postcode{75005}, \country{France}}}


\abstract{
We introduce a new class of large structured random matrices characterized by four fundamental properties which we discuss. We prove that this class is stable under matrix-valued and pointwise non-linear operations. We then formulate an efficient method, based on an extremization problem, for computing the spectrum of subblocks of such large structured random matrices. We present different proofs -- combinatorial or algebraic -- of the validity of this method, which all have some connection with free probability. We illustrate this method with well known examples of unstructured matrices, including Haar randomly rotated matrices, as well as with the example of structured random matrices arising in the quantum symmetric simple exclusion process. 
}
%We present a new efficient method, based on an extremization problem, for computing the spectrum of subblocks of large structured random matrices. This method applies to ensembles of matrices satisfying three fundamental properties which we discuss. We present different proofs -- combinatorial or algebraic -- of the validity of this method, which all have some connection with free probability. We illustrate this method with well known examples of unstructured matrices, including Haar randomly rotated matrices, as well as with the example of structured random matrices arising in the quantum symmetric simple exclusion process. 

\maketitle
\tableofcontents

%\section*{Editorial comments}

\section{Introduction and general statements}\label{sec1}

The theory of large random matrices has a huge domain of applications ranging from chaotic systems to complex systems to random geometry to machine learning \cite{Mehta2004Random,Bouchaud2020First,CouilletLiao2022}. Given a large random matrix $M$ one might not only be interested in its spectrum but also in the spectrum of its subblocks (or submatrices). Moreover, the matrix $M$ might have some "structure", in the sense that joint moments of its entries can depend on the location of these entries inside the matrix. In other words, a structured matrix is, in law, not invariant under permutations of its entries \cite{vanHandel2017Structured} -- contrary to well-known matrix ensembles such as say Wigner matrices.

Finding the spectrum of structured matrices and their subblocks is a problem that can occur in many situations, e.g. in the study of random band matrices \cite{Bourgade2018Random}. Our main original motivation, however, comes from the problem of calculating the entanglement entropy of some many-body quantum systems that are subjected to noise. In this case the system density matrix $\rho$ is a large random matrix, and to calculate the entanglement between a subregion $I$ and the rest of the system $I^c$ requires knowing the so-called reduced density matrix $\rho_I=\Tr_{I^c}(\rho)$. More precisely, we encountered this problem in studying an one dimensional chain of noisy free fermions named the "Quantum Symmetric Simple Exclusion Process" (QSSEP) \cite{Bernard2019Open,Hruza2023Coherent,Bernard2023Exact}. Here, the quadratic (but noisy) Hamiltonian ensures that all properties of the system can be expressed in terms of the two point function $M_{ij}:=\Tr(\rho\, c_i^\dagger c_j)$ where $c_i^\dagger$ is a fermionic creation operator on site $i$. Since the dynamics is noisy, $M$ is a large random matrix and for the entanglement entropy of a region $I$ we need to find the spectrum of any of its subblocks $M_I=(M_{ij})_{i,j\in I}$. The main physical output of the exact computation of \cite{Bernard2023Exact} is that the mutual information in the driven out-of-equilibrium QSSEP fulfills a volume law\footnote{That is, the mutual information between extensive sub-intervals scales proportionally to the volume.}, in contrast with equilibrium systems for which the mutual information is sub-leading in the volume.

Despite this specific motivation, the aim of this paper is to extract from these studies of noisy many-body quantum systems the random matrix related results and to make them  available to a larger audience interested in random matrix theory.
They apply to a large class of ensembles of structured random matrices characterized by specifying the large size limit of so-called "loop expectation values". Loop expectation values are expectation values of the product of entries of random matrices whose indices follow a cyclic order  (see below). It has recently been recognized that these loop expectation values play a peculiar role in abstract random (structured) matrix theory \cite{NICA2002227,Speicher2023NonLinear}, but also in physical contexts \cite{Foini2019Eigenstate,Bernard2019Open,Hruza2023Coherent,Pappalardi2022ETH,Biane2021Combinatorics} or in connection with machine learning \cite{CouilletLiao2022,Speicher2023RMML}. 
Specifying the ensemble of random matrices through the loop expectation values allows to make the connection with the combinatorics of partitions or with free probability transparent and efficient. An echo of this connection is a formula for the moments of the random matrix, see eq.\eqref{eq:n-points-pi} below, as a sum over non crossing partitions and their Kreweras duals. This formula was proved in \cite[sec.\ II.B]{Hruza2023Coherent}. This class and its characterization in terms of loop expectation values is new, to the best of our knowledge. We discuss a few of the properties that this class of ensembles enjoys. In particular, it is stable under non-linear operations as we proved below in Propositions \ref{prop:matrix-poly} and \ref{prop:point-wise}.

The main theorem \ref{thrm:F} stated below reduces the computation of the spectrum of random matrices in such ensembles to a variational problem, which may be viewed as some variante of the Legendre transform, see eq.\eqref{eq:action}, or as a local version of the known R-transform in free probability, see eq.\eqref{eq:local-R-transform}.  Its proof is based on the combinatorial formula \eqref{eq:moment_free_cumulants}. Some elements of those proofs have been evasively formulated in our previous paper \cite{Bernard2023Exact}. We nevertheless believe that it is useful to present synthetically these results and proofs in a separate publication devoted to random matrices -- and not keep them hidden in articles devoted to quantum physics. We also believe that presenting different proofs may be useful depending on the background of the  readers. We complement the two combinatorial proofs by describing a third proof using free probability techniques, notably operator valued free probability or free amalgamation. One the one hand, this makes the connection of the combinatorial proofs with free probability explicit and, on the other hand, it makes concrete the applications of free amalgamation in the present context of structured random matrices. 

The paper is organized as follows. The rest of this section \ref{sec1} summarizes the main properties of this class of random matrices, including a formulation and a discussion of the defining axioms, a description of its stability under non-linear operations, and finally the constructive theorem about the spectrum of sub-blocks. Section \ref{sec2} contains the proofs of the statements made in Section \ref{sec1}. The variational method for computing spectra of sub-matrices is illustrated in Section \ref{sec3} with known examples of rotation invariant matrices and with a new application to the quantum symmetric simple exclusion process.

Before starting, let us fix a point of notation. For random variables $X_1,\cdots,X_n$, we shall denote their moments by $\E[X_1\cdots X_n]$ and their cumulants by $C_n[X_1,\cdots,X_n]$.

\subsection{Axioms and cyclic cumulants} 

Let us start by specifying the random matrix ensembles we shall deal with.
To any monomials in the matrix elements $M_{i_1j_1}\cdots M_{i_nj_n}$ we may associate a graph by assigning to each matrix elements $M_{ij}$ an oriented edge from vertex $i$ to vertex $j$. For example
\begin{align*}
	M_{i_1 i_2}|M_{i_2i_3}|^2M_{i_2i_3}M_{i_3i_1}M_{i_1 i_4}M_{i_4 i_1}M_{i_5 i_5} = \raisebox{-0.5\height}{\includegraphics{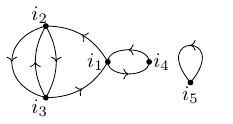}}.
\end{align*}
To a monomial $M_{i_1 i_2}M_{i_2i_3}\cdots M_{i_n i_1}$ with cyclic indices one associates a loop. We can formulate the axioms defining these ensembles either in terms of expectation values of monomials or of graphs.

\paragraph{Random Matrix Ensemble Axioms.} \label{axioms}
We shall consider ensembles of random hermitian matrices $M=M^\dag$ with measure $\mathbb E$ that satisfy, in the large $N$ limit, the three following defining properties~:
\begin{enumerate}
	\item[(i)] Local $U(1)$-invariance, meaning that, in distribution, $M_{ij}\stackrel{d}{=}e^{-i\theta_i}M_{ij}e^{i\theta_j}$, for any phases $\theta_i$, $\theta_j$;
	\item[(ii)] Expectation values of loops without repeated indices scale as $N^{1-n}$, meaning that $\mathbb E [M_{i_1 i_2}M_{i_2i_3}\cdots M_{i_n i_1}]=\mathcal O(N^{1-n})$, for all $i_k$ distinct;
	\item[(iii)] Scaled cumulants of loops $N^{n-1}C_n[M_{i_1 i_2},M_{i_2i_3},\cdots ,M_{i_n i_1}]$ are continuous functions in $x_k=i_k/N$ at coinciding indices, in the large $N$ limit;
	\item[(iv)] Expectation values of disconnected graphs factorize at leading order, meaning that $\mathbb E[M_{i_1 i_2}\cdots M_{i_{m} i_1}\,M_{j_{1} j_{2}}\cdots M_{j_n j_{1}}]=\mathbb E [M_{i_1 i_2}\cdots M_{i_{m} i_1}]\,\mathbb E [M_{j_{1} j_{2}}\cdots M_{j_{n} j_{1}}](1+\mathcal O(N^{-1}))$, where the sets $\{i_1,\cdots,i_n\}$ and $\{j_1,\cdots,j_m\}$ are disjoint, but indices within the sets may coincide. (Note that in \cite{bernard2025additionstructuredrandommatrices} we have proposed a stronger version of axiom (iv) in terms of cumulants  for which all claims made in this paper continue to hold true.)
	
\end{enumerate}

Axiom (i) is clearly weaker than the $U(N)$ invariance in law of usual ensembles of unstructured hermitian random matrices, such as Gaussian Wigner matrices or Haar randomly rotated matrices. This axiom restricts this symmetry to the subgroup $U(1)^N\subset U(N)$, and it is thus compatible with the existence of structure in the matrices. Therefore, the expectation value of a graph is non-vanishing, only if the graph is Eulerian (with identical number of on-going and out-going edges at each vertex).

Axiom (ii) can alternatively be formulated in terms of the cumulants by demanding that $C_n [M_{i_1 i_2},M_{i_2i_3},\cdots ,M_{i_n i_1}]=\mathcal O(N^{1-n})$, since the cumulants and the expectation values of loops coincide for distinct indices. We shall call them \textit{cyclic cumulants} and \textit{loop expectation values}, respectively.
So equivalently, this axiom assumes the finiteness of the scaled cyclic cumulants (with $x_k=i_k/N\in[0,1]$)~:
\begin{equation} \label{eq:def_g}
	g_n(x_{1},\cdots,x_{n}):=\lim_{N\to\infty}N^{n-1}C_n[M_{i_{1}i_{2}}M_{i_{2}i_{3}}\cdots M_{i_{n}i_{1}}] ~.
\end{equation}
For reasons explained below, we shall call the $g_n$ the local free cumulants. Their integrated version (over $x_1,\cdots,x_n$) are not the free cumulants of $M$, but the $g_n$ are related to the operator-valued free cumulants of $M$, as we explain in section \ref{subsec:proof_op-val}.

By Axiom (iii), these local free cumulants are continuous when two or more arguments coincide. For instance $N\,C_2[M_{ii},M_{ii}]=\lim_{y\to x}g_2(x,y)(1+\mathcal O(N^{-1}))$. This property ensures that the expectation value of pinched loops, such as $\E [M_{ii}M_{ii}]$, factorizes at leading order, for instance $\E [M_{ii}M_{ii}]=\E [M_{ii}]^2(1+\mathcal O(N^{-1}))$. More generally, any connected graph can be expressed to leading order in $1/N$ in terms of products of the expectation value of its internal loops. For instance, $\mathbb{E}[M_{12}M_{21}M_{12}M_{21}M_{12}M_{21}]=6\,\mathbb{E}[M_{12}M_{32}]^{3}(1+\mathcal O(N^{-1}))$. This axioms also yields partial information on sub-leading contributions to pinched loop expectation values.

Axiom (iv) and (iii) ensures that the expectation values of any graph can be evaluated in terms of the loop expectation values since any connected graph can be obtained as the limit of a loop with coinciding points. Thus, the only information on the random matrix ensemble we require are the local free cumulants \eqref{eq:def_g}.

Of course, not any sequence of numbers can define cumulants. So we cannot construct an ensemble of random matrices by an arbitrary choice of $g_n$ subjected to the four axioms (i)-(iv). Rather one has to start from a known random matrix ensemble and check that it satisfies (i)-(iv). It remains an open question (at least to us) to list a set of conditions on those $g_n$ to specify a proper well-defined measure on random matrices\footnote{In QSSEP, the local free cumulants $g_n$ can related to the free cumulants of an auxiliary measure. Namely, in this case the series of function $\varphi_n$ defined by $\varphi_n(x_1,\cdots,x_n):=\sum_{\pi\in NC(n)} g_\pi(\vec x)$ are the moments of variables $\mathbb{I}_x$ w.r.t. to an appropriate measure. It remains an open question to know whether this mapping to an auxiliary measure holds true for any ensemble in the class we consider.}.

Some well known matrix ensembles that satisfy these properties are Wigner matrices and matrices rotated by Haar random unitaries (see subsection \ref{subsec:wigner} and \ref{subsec:Haar}), for which the functions $g_n$ are all constant, implying that these ensembles are "structureless". In particular, for $M = UDU^\dagger$ with $D$ a diagonal matrix and $U$ a Haar random unitary, $g_n(x_1,\cdots,x_n)=\kappa_n$ with $\kappa_n$ the free cumulants of the spectral measure of $D$ (see e.g. Thrm. 7.5 in \cite{Speicher2019Lecture} or Appendix \ref{app:free_cum_haar}). For structured ensembles, where the functions $g_n$ are no longer constant, this observation suggests that we could call $g_n$ the "local free cumulants" of $M$.

These ensembles of random matrices are very much related to the framework of free probability theory. 
A first evidence comes from the fact that loop expectation values of $M$ can be decomposed as non-crossing partitions and their Kreweras duals, see \cite[sec.\ II.B]{Hruza2023Coherent}, (with $\vec{x}=(x_{1},\cdots,x_{n})$ and $x_k=i_k/N$ fixed in the large $N$ limit),
\begin{equation} \label{eq:n-points-pi}
\lim_{N\to\infty} N^{n-1}\mathbb E[ M_{i_1i_2}M_{i_2i_3}\cdots M_{i_ni_1}] =\sum_{\pi\in NC(n)}g_{\pi^*}(\vec x)\,\delta_\pi(\vec x) ~,
\end{equation}
where  $NC(n)$ denotes the set of non-crossing partitions of order $n$, and $g_\pi(\vec x):=\prod_{p\in\pi}g_{|p|}(\vec x_p)$ with $\vec x_p=(x_i)_{i\in p}$ the collection of variables $x_i$ belonging to the part $p$ of the partition $\pi$, and $|p|$ the number of elements in this part. By $\delta_\pi(\vec x)$ we denote a product of delta functions $\delta(x_i-x_j)$ that equate all $x_i,x_j$ with $i$ and $j$ in the same part $p\in \pi$. And $\pi^*$ is the Kreweras complement of $\pi$ (see Section \ref{sec:subblocks} for an example and \cite{Speicher2019Lecture} for the definition of the Kreweras dual). 
Of course, evaluating the cumulants of $M$ from this expression yields back \eqref{eq:def_g}.

To prevent a confusion, note that $\tilde \kappa_\pi:=\int g_{\pi^*}(\vec x)\,\delta_\pi(\vec x) d\vec x$ are not the free cumulants of $M$, because they fail to be multiplicative, i.e. $\tilde \kappa_\pi \tilde\kappa_\sigma\neq\tilde\kappa_{\pi\cup\sigma}$ with $\pi\cup\sigma$ the union of parts of $\pi$ and $\sigma$. The reason for this is the contraction with the delta function\footnote{Formally, the free cumulants of $M_h$ are defined as a multiplicative family $\kappa_\pi=\prod_{b\in\pi} \kappa_{|b|}$ with $\kappa_n:=\kappa_{1_n}$ satisfying $\phi_n[h] = \sum_{\pi\in NC(n)}\kappa_\pi$ and they can be related to the $\tilde \kappa_\pi$ by Moebius inversion, $\kappa_n = \sum_{\pi\in NC(n)} \mu(\pi,1_n) \prod_{b\in\pi}\sum_{\sigma\in NC(|b|)}\tilde \kappa_\sigma$.}.

According to the above axioms, expectation values of any monomials in the matrix elements of $M$ are expressible in terms of local free cumulants \eqref{eq:def_g}. This property yields an explicit formula for the cumulant generating function\footnote{Note that formally in the scaling limit, $\mathbb{E}[e^{N{\ntr}(MQ)}]=\mathbb{E}[e^{N\!\int\!m(y,x)q(x,y)dxdy}]$ with $m(x,y)=M_{ij}$ in the large $N$ limit.} $\log\mathbb{E}[e^{N\,{\ntr}(MQ)}]$, with $\ntr(\cdot)=(1/N)\mathrm{tr}(\cdot)$ the normalized trace, for some test matrix $Q$. We choose $Q$ to have a proper large $N$ limit.
% so that $q(x,y)=\lim_{N\to\infty} Q_{ij}$ exits, for $i=[x/N],\ j=[y/N]$.
\medskip

\begin{proposition} 
Let $Q_{ij}=q(\frac{i}{N},\frac{j}{N})$ for some regular enough function $q(x,y)$. Then, 
\[
%\mathbb{E}[e^{N \mathrm{tr}(MQ)}] \stackrel{N\to\infty}{\asymp}e^{N W[q]}
\mathbb{E}[e^{N {\ntr}(MQ)}]\, {\asymp}_{N\to\infty}\, e^{N W[q]}
\]
with the cumulant generating function $W[q]$  given by (at least as a formal power series)
%, and with $\vec{x}=(x_{1},\cdots,x_{n})$) 
\begin{equation} \label{eq:W-cumulant}
W[q] =\sum_{n\geq 1}\frac{1}{n}\int_0^1 \! g_{n}(x_{1},\cdots,x_{n})\,q(x_{1},x_{n})\cdots q(x_{3},x_{2})q(x_{2},x_{1})\,d\vec{x} ~.
\end{equation}
\end{proposition}
%\medskip

See the proof in Section \ref{sec:generating-function}.
In the case of unstructured matrices, with the loop expectation values $g_n$ independent of the positions, i.e. $g_{n}(\vec{x})=\kappa_n$, we simply have $W[q]=\sum_{n}\frac{1}{n}\kappa_{n}\mathrm{tr}(\hat Q^{n})$, with $\hat Q=Q/N$. For Haar randomly rotated matrices, $\kappa_n$ is the $n$-th free cumulants of the spectral measure of the matrix $D$.

In the special case with $q(x,y)=h^{\frac{1}{2}}(x)h^{\frac{1}{2}}(y)$, for some function $h$, we are dealing with the matrix element $v_h^TMv_h$, with $v_h$ the vector with coordinates $h_i^\frac{1}{2}$. Equation \eqref{eq:W-cumulant} then yields a formula for its cumulant generating function $F_0[h]:=\lim_{N\to\infty}N^{-1}\log\mathbb{E}[e^{v_h^TMv_h}]$,
\begin{equation}
F_0[h] =\sum_{n\geq 1}\frac{1}{n}\int_0^1 \! g_{n}(x_{1},\cdots,x_{n})\,h(x_{1})\cdots h(x_{n})\,d\vec{x} ~.
\end{equation}
Note that, contrary to $W[q]$ in \eqref{eq:W-cumulant}, $F_0[h]$ encodes information only on the permutation invariant part of the local free cumulants. This function is going to play an important role for determining the spectrum of $M_h$.

\subsection{Closure under non-linear operations}
The set of random matrix ensembles satisfying the above axioms is stable under some non-linear operations. Combination of these operations provides in particular a way to generate ensembles of structured random matrices from unstructured ones.
The first statement concerns matrix valued polynomial operations.
\medskip

\begin{proposition} \label{prop:matrix-poly}
The axioms (i)-(iv) close under matrix valued polynomials.
That is: if $M$ satisfy these axioms, so does $P(M)$ for any polynomial $P(M)=\sum_{m\geq 0} a_mM^m$.\\
If $g_n$ are the local free cumulants of $M$, the local free cumulants of $P(M)$ (defined as $g_n^P(x_1,\cdots,x_n):=\lim_{N\to\infty} N^{n-1} C_n[P(M)_{i_1i_2},\cdots,P(M)_{i_ni_1}]$) are given as
\begin{align} \label{eq:g-polynomials}
g_{n}^{P}(x_{1},\cdots,x_{n})
 & =\sum_{m_{1},\cdots,m_{n}}a_{m_{1}}\cdots a_{m_{n}} \hskip -0.2 truecm \sum_{\pi\in NC(m)\atop\text{s.t. }\pi\vee\Gamma=1_{m}}  \int_0^1 g_{\pi}(x_{\cdot},\vec{y}_{\cdot})\delta_{\pi^{*}}(x_{\cdot},\vec{y}_{\cdot})d\vec{y}_{\cdot} ~,
\end{align}
where in each term of the sum,  $m=m_{1}+\cdots+m_{n}$ and $\Gamma=\gamma_1\cup\cdots\cup\gamma_n$ is a partition of $\{1,\cdots,m\}$ into intervals of length $|\gamma_i| = m_i$. Furthermore, $\pi\vee \Gamma$ denotes the common maximum of $\pi$ and $\Gamma$ in the partially ordered set of non-crossing partitions. We used $\int g_{\pi}(x_{\cdot},\vec{y}_{\cdot})\delta_{\pi^{*}}(x_{\cdot},\vec{y}_{\cdot})d\vec{y}_{\cdot}$ as a shorthand notation for $\int g_{\pi}(x_{1},\vec{y}_{1},\cdots,x_{n},\vec{y}_{n})\delta_{\pi^{*}}(x_{1},\vec{y}_{1},\cdots,x_{n},\vec{y}_{n})d\vec{y}_{1}\cdots d\vec{y}_{n}$ with integration variables $\vec y_i = (y_1^{(i)},\cdots,y^{(i)}_{m_i})$.
\end{proposition}
\medskip

The proof is given in Section \ref{sec:operation}. Let us here illustrate this statement on a simple example with $P(M)=M^2$. We compute the first few loop expectation values of $P(M)$ when indices $x_k=i_k/N$ are distinct. They are thus equal to the local free cumulants $g^P_n$,
\begin{align*}
g_1^P(x) &= \sum_k \mathbb{E}[M_{ik}M_{ki}] = g_1(x)^2+\int_0^1\!g_2(x,y)dy ,  \\
g_2^P(x_1,x_2)&= N\sum_{jk}\mathbb{E}[M_{i_1j}M_{ji_2}M_{i_2k}M_{ki_1}] \\
& = \int_0^1\!g_4(x_1,y_1,x_2,y_2)dy_1dy_2 + \int_0^1\! g_1(x_1)g_3(x_1,x_2,y)dy  \\
&~~~ + \int_0^1\! g_3(x_1,y,x_2,)g_1(x_2)dy  + g_1(x_1)g_2(x_1,x_2)g_1(x_2) .
\end{align*}
The key point is that the scaling of the cumulants with $N$ is compensated by the summation over the intermediate indices in the matrix products.

The second statement concerns point-wise non-linear operations, which are similar to those involved in large neural networks  \cite{CouilletLiao2022,NIPS2017_0f3d014e,10.1214/21-EJP699,Speicher2023RMML}. The statement requires $M$ to be centred, $\E[M_{ii}]=0$. This can be assumed without loss of generality, since the axioms are preserved by shifting $M$ by a diagonal matrix, $M_{ij}\to M_{ij} - \delta_{ij} a_i$.
\medskip

\begin{proposition} \label{prop:point-wise}
Given $M$ with $\E[M_{ii}]=0$, let us define a new matrix $Y_{ij}= M_{ij}\, f_{ij}(N|M_{ij}|^2)$ where $f_{ij}(\cdot)$ is a non-linear function (defined as power series) that can be different for each entry, with a finite large $N$ limit $f_{xy}(\cdot)$. Then, if $M$ satisfies the axioms (i)-(iv) so does $Y$. \\
For $f_{xy}(u)=\sum_{m\geq 0} f_m u^m$, let us define $D_{f_{xy}}[u]:=\sum_{m\geq 0} (m+1)!f_m u^m$. Then, if $g_n$ are the local free cumulants of $M$, those of $Y$ are 
\begin{equation} \label{eq:gY}
g_n^Y(x_1,\cdots,x_n) = g_n(x_1,\cdots,x_n)\, D_{f_{x_1,x_2}}[g_2(x_1,x_2)]\cdots D_{f_{x_n,x_1}}[g_2(x_n,x_1)] ~.
\end{equation}
\end{proposition}
%\medskip

The proof is given in Section \ref{sec:operation} and we illustrate it here on a simple example. Let $f_{ij}=\text{id}$, so $Y_{ij}= N\, M_{ij} |M_{ij}|^2$, and let us compute the loop expectation values $\mathbb{E}[Y_{i_1i_2}\cdots Y_{i_ni_1}]$, with all $i_k$ distinct. Due to $U(1)$ invariance, the leading contributions to this expectation value come from partitioning the product of matrices $M$ into one group of size $n$ made of the cyclic product  $M_{i_1i_2}\cdots M_{i_ni_1}$ and $n$ other groups each of size $2$ made of $|M_{i_ki_{k+1}}|^2$. There are $2^n$ manners of partitioning this way. Thus, to leading order in $1/N$, (with $x_k=i_k/N$ and all $i_k$ distinct),
\begin{align*}
\mathbb{E}[Y_{i_1i_2}\cdots Y_{i_ni_1}] &= 2^n\, N^n\, \mathbb{E}[M_{i_1i_2}\cdots M_{i_ni_1}]\, \mathbb{E}[|M_{i_1i_2}|^2]\cdots \mathbb{E}[|M_{i_ni_1}|^2]\, (1+\mathcal O(N^{-1}))\\
&= N^{1-n}\, 2^n\, g_n(x_1,\cdots,x_n)\, g_2(x_1,x_2)\cdots g_2(x_n,x_1)\, (1+\mathcal O(N^{-1})).
\end{align*}
In particular, the $n$-point functions of the $Y$'s have the appropriate scaling as function of $N$. 
It is also simple to verify that $Y$ satisfies the other axioms.

Note that formula \eqref{eq:gY} amounts to dress the edges of the loop associated to the cyclic product $Y_{i_1i_2}\cdots Y_{i_ni_1}$ by the factors $D_{f_{x_k,x_{k+1}}}[g_2(x_k,x_{k+1})]$. This simple observation has the direct consequence that the cumulant generating function of $Y$ can be obtained from that of $M$ via a simple dressing. Namely, if $W_M$ is the cumulant generating function of $M$, as in \eqref{eq:W-cumulant}, that of $Y$ is
\begin{equation}
W_Y[q] = W_M[ q^{D_f}] ~,
\end{equation}
where the dressed kernel is $q^{D_f}(y,x) = q(y,x)\,D_{f_{x,y}}[g_2(x,y)]$.

To explore the space of possible local free cumulants, one can -- as stated above -- combine the two above propositions and generate new ensembles of structured random matrices from unstructured once. One may for instance start with perturbed Wigner random matrices with measure $dM\, e^{-N \tr V(M)}$ for some potential $V(M)$, in the large $N$ limit. Then, iterating (a) entry-dependent and (b) matrix-valued non-linear operations, one obtains new ensembles of structured random matrices satisfying the axioms (i)-(iv).

 \subsection{Spectrum of subblocks}

To handle the case of an arbitrary number of subblocks we consider the slightly more general aim of finding the spectrum of $M_h:=h^{\frac{1}{2}}Mh^{\frac{1}{2}}$, with $h$ a diagonal matrix. Choosing $h(x)=1_{x\in I}$ (here $h_{ii}=h(i/N)$) to be the indicator function on some interval $I\subset[0,1]$, one recovers the case of subblocks $M_I\subset M$. All the spectral information about $M_h$ is contained in the
generating function,
\begin{align}\label{eq:def_F}
	F[h](z):=\mathbb{E}\,\ntr\log(z-M_h),
\end{align}
where $\ntr=\tr/N$ is the normalized $N$-dimensional trace. 
This function can be viewed as a (formal) power series in $1/z$, whose coefficients are expectations of traces of powers of $M_h$. Statements about the domain of convergence of this series can be made if extra global information about the spectrum is available, say about its compactness. The theorem below is formulated with $F[h](z)$ viewed as  power series in $1/z$ (and we use extra analytic inputs in the illustrative examples).

Our main result is~:
\medskip

\begin{theorem}\label{thrm:F}
$F[h](z)$ is determined by the variational principle
\begin{equation}\label{eq:action}
	F[h](z)=\underset{a_z,b_z}{\mathrm{extremum}}\left[\int_0^1\left[\log(z-h(x)b_z(x))+a_z(x)b_z(x)\right]dx-F_{0}[a_z]\right]
\end{equation}
where the information about local free cumulants $g_n$, specific to the random matrix ensemble, is contained in (with $\vec x=(x_1,\cdots,x_n)$) 
\begin{equation}
	F_0[p]:=\sum_{n\ge1}\frac{1}{n}\int_0^1 (\prod_{k=1}^ndx_kp(x_{k}))\, g_n(\vec x).
\end{equation}
\end{theorem}

Note that $F_0[p]$ contains less information than the local free cumulants, since it depends only on a symmetrized version of the family $\{g_n\}_n$. Nevertheless, in the large $N$ limit, it represents the minimal amount of information about the measure $\mathbb{E}$ that is necessary for the spectrum.

To obtain the spectrum of $M_h$ one takes the derivative $\partial_z F[h](z)=: G[h](z)$ which is the resolvent
\begin{equation*}  %\label{eq:resolvent}
G[h](z) = \mathbb E\,\ntr {(z-M_h)^{-1}} . %\overset{!}{=}\int_0^1 \!\frac{dx}{z-h(x)b_z(x)}
\end{equation*}
From Eq.\eqref{eq:action}, we get
 \begin{equation}\label{eq:resolvent}
G[h](z)=\int_0^1 \!\frac{dx}{z-h(x)b_z(x)} ,
\end{equation}
with $b_z$ solution of the extremization conditions,
\begin{equation} \label{eq:a-b}
	a_z(x)=\frac{h(x)}{z-h(x)b_z(x)},\quad b_z(x)= R_0[a_z](x) ,
\end{equation}
where
\begin{equation}\label{eq:R_0}
	R_0[a_z](x):= \frac{\delta F_0[a_z]}{\delta a_z(x)}.
\end{equation} 

In the special case where $h(x)=1_{x\in I}$ is the indicator function on an interval $I$ (or on unions of intervals) of length $\ell_I$, we recover the normalized spectral density $\sigma_I$ of the subblock $M_I$ from its resolvent 
\begin{equation}\label{eq:resolvent_I}
	G_I(z):=\frac{1}{\ell_I}\int_I \frac{d x}{z-b_z(x)}=\int \frac{d\sigma_I(\lambda)}{z-\lambda}
\end{equation}
as $G_I(\lambda-i\epsilon) -  G_I(\lambda+i\epsilon)= 2i\pi\sigma_I(\lambda)$.
Writing the total resolvent (including the pole at the origin)
\[
G^\mathrm{tot}_I(z):=G[1_{x\in I}](z)= \frac{1-\ell_I}{z} + \ell_I \int \frac{d\sigma_I(\lambda)}{z-\lambda},
\]
we can relate the total spectral measure of $M_h$ (including the zero-eigenvalues) to that of a subblock $M_I\subset M$ by \begin{equation}\label{eq:total_spectrum}
	d\sigma^\mathrm{tot}_I(\lambda)=(1-\ell_I)\delta(\lambda)d\lambda+ \ell_I\, d\sigma_I(\lambda).
\end{equation}

This result is of course very much related to the framework of free probability theory, as illustrated by the following rewriting of Eq.\eqref{eq:a-b},
\begin{equation} \label{eq:local-R-transform}
zh(x)^{-1}=a_z(x)^{-1}+R_0[a_z](x).
\end{equation}
This resembles a local version of the so-called R-transform of free probability theory -- hence the choice for the name of "local free cumulants" for $g_n$.

Finally, it turns out that our result can also be obtained from the general relation between the R-transform and the resolvent in the framework of operator-valued free probability theory (see section \ref{subsec:proof_op-val}). However, the general form of this relation is quite abstract and some work is necessary to see that it can be applied to the more practical problem of finding the spectrum of subblocks of a class of random matrices satisfying axioms (i)-(iv). This is one of the main contributions of this article -- a second being a direct proof of our result that does not use operator-valued free probability theory (see sections \ref{subsec:proof_trees} and \ref{subsec:proof_kreweras})

Besides the general statement, the application of our method to the QSSEP random matrix ensemble also constitutes a new result (see section \ref{subsec:QSSEP}). In contrast to this, the other applications we present (see sections \ref{subsec:wigner} and \ref{subsec:Haar}) are rather illustrations on how to use our method in some well-known matrix ensembles in order to make contact with known results.
	
Recalling that we parametrize an arbitrary collection of subblocks of $M$ by the action with a diagonal matrix $h^{\frac{1}{2}}$, i.e. $M_h=h^{\frac{1}{2}}M h^{\frac{1}{2}}$, one might wonder if the spectrum of $M_h$ can be obtained (much faster) by free convolution~? Indeed, as long as $M$ and $h$ are free in the sense of (scalar) free probability theory, the spectral measure of $M_h$ can be obtained by free multiplicative convolution from the spectral measures of $M$ and $h$ (see e.g. Lecture 14 in \cite{Nica2006Lecture}). But this hope turns out to be in vain (see section \ref{subsec:no_free_mulit_convolution})~: The spectrum of $M_h$ does not coincide with the spectrum obtained by free multiplicative convolution from $M$ and $h$. In turn, this means that structured random matrices satisfying (i)-(iv) are not free from diagonal deterministic ones -- highlighting the special role of structured random matrices.

However, there are special cases of structureless matrix ensembles that are free from diagonal deterministic matrices, for example if $M$ is a matrix rotated by Haar random unitaries (Theorem 7.5 in \cite{Speicher2019Lecture}). In this case we show in section \ref{subsec:Haar} (rather for illustrating than for original purposes) that the spectrum of $M_h$ can be indeed obtained from free convolution. In the same section we also show, that for a single subblock our result from Theorem \ref{thrm:F} reduces to "free compression". For Haar-randomly rotated matrices, this is a well known result.

Let us also note that we can invert the variational principle~: Given a generating function $F[h](z)$ that satisfies Eq.\eqref{eq:action}, we can retrieve the initial data $F_0$ as the extremum of
\begin{equation}\label{eq:inverted_action}
	F_0[a]=\underset{h,b_z}{\mathrm{extremum}}\left[\int_0^1\left[\log(z-h(x)b_z(x))+a_z(x)b_z(x)\right]dx-F[h](z)\right].
\end{equation}
This is very similar to the Legendre Transformation where the initial function can be retrieved by applying the transformation twice.
Here the inversion works because in extremizing Eq.\eqref{eq:action} we obtain $a=a(h,z)$ and $b=b(h,z)$ as functions of $h$ (and $z$), while  in extremizing Eq.\eqref{eq:inverted_action} we obtain $h=h(a,z)$ and $b=b(a,z)$ as functions of $a$ (and $z$). Through formal power series, the triple $(a,b,h)$ can be inverted which ensures the variational principle for $F_0$ above.

\section{Proofs} \label{sec2}

\subsection{Generating functions}
\label{sec:generating-function}
Here we give the proof of the formula \eqref{eq:W-cumulant} for the cumulant generating function in terms of the local free cumulants.
	
Consider $X=\ntr(MQ)$ so that $NX= \tr(MQ)$. The classical moment-cumulant relation allows us to write
\begin{equation*}
	NW[q]=\sum_{n\ge1} \frac{1}{n!} \sum_{ij}C_n[M_{i_1j_1},\cdots ,M_{i_nj_n}]\, Q_{j_ni_n}\cdots Q_{j_1i_1} ~,
%	NW[q]=\sum_{n\ge1} \frac{1}{n!}\, \E[(\tr(MQ))^n]^c=\sum_{n\ge1} \frac{1}{n!} \sum_{ij}C_n[\hat M_{i_1j_1},\cdots ,\hat M_{i_nj_n}],
\end{equation*}
where we abbreviate $i=(i_1,\cdots,i_n)$ and $j=(j_1,\cdots,j_n)$. In the sum over $i$, one can keep only distinct indices, taking away all the other terms where two or more indices are equal since those terms introduce an error that is sub-leading in $1/N$. Then, due to the $U(1)$ invariance of the measure, one has (with the shortened notation $\hat M_{ij} := M_{ij} Q_{ji}$)
\begin{equation}
	\sum_{i \text{ distinct}} \sum_{j}C_n[\hat M_{i_1j_1},\cdots ,\hat M_{i_nj_n}]= \sum_{i \text{ distinct}} \sum_{\sigma\in S_n}C_n[\hat M_{i_1i_{\sigma(1)}},\cdots ,\hat M_{i_ni_{\sigma(n)}}] ~,
\end{equation}
where $\sigma\in S_n$ is a permutation of $n$ elements. Permutations consisting of a complete cycle such as $\sigma=(1\cdots n)$ produce terms of the form $C_n[M_{i_1i_2},\cdots, M_{i_ni_1}]\sim\Ord(N^{1-n})$, while all other permutations, consisting of more than one cycle, produce sub-leading terms. For example $\sigma=(1)(2\cdots n)$ leads to $C_n[M_{i_1i_1}, M_{i_2 i_3},\cdots ,M_{i_ni_2}]\sim\Ord(N^{-n})$.
Therefore, one keeps the $(n-1)!$ complete cycles which all give the same contribution, since we can order the monomials as $M_{i_1i_{\sigma(1)}}M_{i_{\sigma(1)}i_{\sigma^2(1)}}\cdots$. Thus,
\begin{equation}
	NW[Q] = \sum_{n\ge1} \frac{1}{n} \sum_{i \text{ distinct}} C_n[M_{i_1i_2},\cdots ,M_{i_ni_1}]\, Q_{i_1i_n}\cdots Q_{i_2i_1} ~.
\end{equation}
Recall now that $C_n[M_{i_1i_2},\cdots ,M_{i_ni_1}]=N^{1-n}g_n(x_1,\cdots,x_n)$ with $x_k=i_k/N$.
Allowing terms where indices $i$ are equal will again only make a sub-leading error, by continuity of the local free cumulants. Therefore, one can replace the sum by an integral, using the scaling of $Q_{ij}=q({i}/{N},{j}/{N})$, which leads to \eqref{eq:W-cumulant}.

\subsection{Non-linear operations}
\label{sec:operation}

We here give the proofs that the set of random matrix ensembles satisfying the above axioms is stable under non-linear operations.

\subsubsection{Polynomial operations}

Let us prove that the axioms (i)-(iv) close under matrix-valued polynomial operations. 
Let us write an arbitrary polynomial in $M$ as $P(M)=\sum_{m>0}a_{m}M^{m}$.
%We would like to show that this satisfies axioms (1), (2') and (3').
Due to multi-liniarity of the cumulants we have
\[
C_{n}[P(M)_{i_{1}i_{2}},\cdots,P(M)_{i_{n}i_{1}}]=\sum_{m_{1},\cdots,m_{n}}a_{m_{1}}\cdots a_{m_{n}}\,C_{n}[(M^{m_{1}})_{i_{1}i_{2}},\cdots,(M^{m_{n}})_{i_{n}i_{1}}] ~.
\]
Denoting the indices that are summed over schematically by $i'$, we have to evaluate (note the positions of the commas which specify the variables for which we are evaluating the cumulants)
\begin{equation} \label{eq:product-M}
C_{n}[(M^{m_{1}})_{i_{1}i_{2}},\cdots,(M^{m_{n}})_{i_{n}i_{1}}]
=\sum_{i'} C_{n}[\underbrace{M_{i_{1}i'}\cdots M_{i'i_{2}}}_{m_{1}\text{ times}},\cdots,\underbrace{M_{i_{n}i'}\cdots M_{i'i_{1}}}_{m_{n}\text{ times}}] ~.
\end{equation}
%Again, we need to use the combinatorial formula, that relates cumulants of products of variables to cumulants of the individual variables.

We now use the formula between cumulants of variables grouped into products and cumulants of the individual variables, see e.g. \cite{Leonov1959OnAM}, which we recall.
Denote $k=k_{1}+\cdots+k_{n}$ and $\Gamma=\gamma_{1}\cup\cdots\cup\gamma_{n}\in P(k)$ a partition of $\{1,\cdots,k\}$ into intervals $\gamma_{1},\cdots,\gamma_{n}$ of length $k_{1},\cdots,k_{n}$. For $\pi\in P(k)$, denote $\pi\vee\Gamma$ the maximum of the two partitions, then
\begin{equation} \label{eq:cumulants-product}
C_{n}[X_{1}^{k_{1}},\cdots,X_{n}^{k_{n}}]=\sum_{\pi\in P(k):\,\pi\vee\Gamma=1_{k}} \!\! C_{\pi}[\underbrace{X_{1},\cdots,X_{1}}_{k_{1}\text{ times}},\cdots,\underbrace{X_{n},\cdots,X_{n}}_{k_{n}\text{ times}}] ~.
\end{equation}
In other words, one has to sum over all partitions $\pi$ (crossing and non-crossing) that connect the intervals $\gamma_{1},\cdots,\gamma_{n}$ to become the complete interval $1_{k}$. 

In the present case \eqref{eq:product-M}, each product has $m_{i}$ variables and there are $m=m_{1}+\cdots+m_{n}$ individual variables in total. With $\Gamma=\gamma_{1}\cup\cdots\cup\gamma_{n}$ a partition of $\{1,\cdots,m\}$ into intervals $\gamma_{i}$ of size $m_{i}$, we have
\begin{align*}
&C_{n}[M_{i_{1}i'}\cdots M_{i'i_{2}},\cdots,M_{i_{n}i'}\cdots M_{i'i_{1}}]\\
& =\sum_{\pi\in P(m):\,\pi\vee\Gamma=1_{m}} \!\! C_{\pi}[M_{i_{1}i'},\cdots,M_{i'i_{2}},\cdots,M_{i_{n}i'},\cdots,M_{i'i_{1}}] ~.
\end{align*}
From the $U(1)$ invariance plus the scaling properties of the cumulants, we know that, for any given partition $\pi\in P(m)$, only non-crossing partitions contribute to the above cumulants at leading order in $1/N$ when summed over all the indices $i'$, compare to eq.~\eqref{eq:n-points-pi}. The constraint of $U(1)$ invariance can be enforced by multiplying be appropriated Kronecker deltas. That is,
\[
C_{\pi}[M_{i_{1}i'},\cdots,M_{i'i_{2}},\cdots,M_{i_{n}i'},\cdots,M_{i'i_{1}}]\,\delta_{\pi^{*}}(i_{1},i',\cdots,i_{n},i')\sim N^{1-m} ~.
\]
Taking the sum $\sum_{i'}$ compensates a factor of $N^{m-n}$. Therefore, the cumulant at order $n$ of a polynomial of $M$ scales as $N^{1-n}$,
\begin{align*}
 & C_{n}[(M^{m_{1}})_{i_{1}i_{2}},\cdots,(M^{m_{n}})_{i_{n}i_{1}}]\\
 & =N^{1-n}\sum_{\pi\in NC(m)\atop \pi\vee\Gamma=1_{m}}\int d\vec{y}_{1}\cdots d\vec{y}_{n}\, g_{\pi}(x_{1},\vec{y}_{1},\cdots,x_{n},\vec{y}_{n})\,\delta_{\pi^{*}}(x_{1},\vec{y}_{1},\cdots,x_{n},\vec{y}_{n}) ~,
\end{align*}
where $\vec{y}_{i}=(y^{(i)}_{1},\cdots,y^{(i)}_{m_{i}-1})$. Putting all  together yields \eqref{eq:g-polynomials}.
This makes it evident that any polynomial $P(M)$ satisfies axiom (i)-(iii). It is easy to check that the axiom (iv) is also fulfilled.

\subsubsection{Point-wise operations}

Again by multi-linearity of the cumulants, it is enough to prove Proposition \ref{prop:point-wise} for monomials $x^{m_k}$ with different
powers $m_{1},\cdots,m_{n}$ for each entry $M_{i_{1}i_{2}},\cdots,M_{i_{n}i_{1}}$, that is for
\[
\Sigma_{m_1,\cdots,m_n}:= N^{m}\, C_{n}[M_{i_{1}i_{2}}|M_{i_{1}i_{2}}|^{2m_{1}},\cdots,M_{i_{n}i_{1}}|M_{i_{n}i_{1}}|^{2m_{n}}] ~.
\]
with $m=m_1+\cdots+m_n$.
Again, we are going to use the formula \eqref{eq:cumulants-product} for the cumulants of product of variables, with $k_{1}=2m_{1}+1,\cdots,k_{n}=2m_{n}+1$. 
Due to the U(1) invariance, the relevant partitions consist of one cyclic block which takes one element from each interval, and all other blocks are
pairs of elements within an interval. Other partitions yield sub-leading contributions. Thus, $\Sigma_{m_1,\cdots,m_n}$ is equal to (up to sub-leading terms in $1/N$)
\begin{align*}
%& ~~~~~~~\Sigma_{m_1,\cdots,m_n}\\
 & = N^{m}\, \lambda_{\{m_{\cdot}\}}\,  C_{n}[M_{i_{1}i_{2}},\cdots,M_{i_{n}i_{1}}]
C_{2}[M_{i_{1}i_{2}},M_{i_{2}i_{1}}]^{m_{1}}\cdots C_{2}[M_{i_{n}i_{1}},M_{i_{1}i_{n}}]^{m_{n}}\\
&=  N^{1-n}\, \lambda_{\{m_{\cdot}\}}\, g_{n}(x_{1},\cdots,x_{n})\,g_{2}(x_{1},x_{2})^{m_{1}}\cdots g_{2}(x_{n},x_{1})^{m_{n}} ~.
%&=  N^{1-n}g_{n}(x_{1},\cdots,x_{n})\,(m_{1}+1)!g_{2}(x_{1},x_{2})^{m_{1}}\cdots(m_{n}+1)!g_{2}(x_{n},x_{1})^{m_{n}}.
\end{align*}
The combinatorial prefactor $\lambda_{\{m_{\cdot}\}}$ codes for the number of ways of partitioning in one cycle of size $n$ and $n$ other blocks of size $2$. It is given by
\[
\lambda_{\{m_{\cdot}\}}=(m_{1}+1)!\cdots(m_{n}+1)! ~.
\]
Indeed, to build the cyclic block, one has a choice between $m_{1}+1$ elements $M_{i_{1}i_{2}}$ in the product $M_{i_{1}i_{2}}|M_{i_{1}i_{2}}|^{2m_{1}}$ (which corresponds to $\gamma_{1}$ in \eqref{eq:cumulants-product}).
Once this element is chosen, there are $m_{1}!$ possibilities to form pairs with the remaining elements while respecting the U(1) invariance. This makes $(m_{1}+1)!$ possible choices for the first product. The other products follow analogously. Multiplying all possibilities gives $\lambda_{\{m_{\cdot}\}}$. Resuming these contributions, using the multi-linearity of the cumulants, yields \eqref{eq:gY}.
This makes it clear that $Y$ satisfies axiom (i)-(iii). It is easy to check that the axiom (iv) is also satisfied.

\subsection{Spectrum of subblocks}
\label{sec:subblocks}

We propose three proofs for Theorem \ref{thrm:F}. The first is a direct proof that uses a bijection between non-crossing partitions and trees. The second is also direct and is based on organizing the summation on partitions according to their cardinal or that of their Kreweras dual. The third relies on operator valued free probability theory and shows that our result is a special case of the relation between operator-valued R-transform and resolvent. Of course, the three proofs have some elements in common.

Let $M_h:=h^{\frac{1}{2}}Mh^{\frac{1}{2}}$ as in \eqref{eq:def_F}. From \eqref{eq:n-points-pi}, the moments of $M_h$ can be expressed as a sum over non-crossing partitions (with $\ntr=\tr/N$)
%(here $NC(n)$ denotes the set of non-crossing partitions of order $n$ and $\vec{x}=(x_{1},\cdots,x_{n})$)
\begin{equation}\label{eq:moment_free_cumulants}
	\phi_n[h]:= \lim_{N\to \infty} \mathbb E[ \ntr(M_h^n)]=\sum_{\pi\in NC(n)}\int \! g_{\pi^*}(\vec x)\,\delta_\pi(\vec x)\, h(x_1)\cdots h(x_n)\,d\vec x
\end{equation}
where the notation is as in \eqref{eq:n-points-pi}.
%where $g_\pi(\vec x):=\prod_{p\in\pi}g_{|p|}(\vec x_p)$ with $\vec x_p=(x_i)_{i\in p}$ the collection of variables $x_i$ belonging to the part $p$ of the partition $\pi$, and $|p|$ the number of elements in this part. By $\delta_\pi(\vec x)$ we denote a product of delta functions $\delta(x_i-x_j)$ that equate all $x_i,x_j$ with $i$ and $j$ in the same part $p\in \pi$. And $\pi^*$ is the Kreweras complement of $\pi$ (see Section \ref{sec:subblocks} for an example and \cite{Speicher2019Lecture} for the definition of the Kreweras dual). 

\subsubsection{Proof using a tree structure}\label{subsec:proof_trees}
Expanding the generating function \eqref{eq:def_F} in terms of the moments $\phi_n[h]$, defined in Eq.\eqref{eq:moment_free_cumulants}, one has
\begin{equation}\label{eq:F_expansion}
	F[h](z)=\log(z)-\sum_{n\ge1} \frac{z^{-n}}{n}\phi_n[h].
\end{equation}
The difficulty in this function lies in organising the sum over non-crossing partitions of any possible integer $n$. To better understand this structure, we note that non-crossing partitions $\pi\in NC(n)$ are in one-to-one correspondence with planar bipartite rooted trees with $n$ edges, if one labels its black and white vertices by the parts of $\pi$ and $\pi^*$. Here is an example for $\pi=\{\{1,3\},\{2\},\{4,5\},\{6\}\}$ (dotted lines) whose Kreweras complement is $\pi^*=\{\{\bar 1, \bar 2\},\{\bar3,\bar5,\bar6\},\{\bar 4\}\}$ (solid lines).
\begin{center}
	\includegraphics[width=0.5 \linewidth]{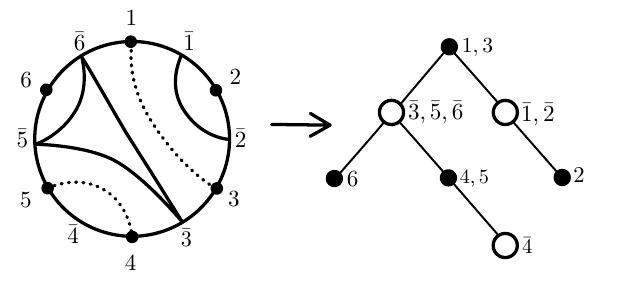}
\end{center}
The parts of $\pi$ are associated with black vertices and parts of $\pi^*$ with white vertices. Two vertices are connected if the corresponding parts of $\pi$ and $\pi^*$ have an element in common (identifying numbers with and without bar, $k\sim\bar k$). The root is (by convention) chosen to be the part $p$ containing $1$.

However, applying this correspondence to Eq.\eqref{eq:moment_free_cumulants} is not directly straightforward, because two partitions $\pi$ and $\pi'$ that are related by a rotation of its elements (in the circle representation) have the same contribution in the sum and thereby complicate the counting of terms. This is due to the integration over $x_1,\cdots, x_n$. If instead, we don't integrate over one of these variables, call it $x$, then $\pi$ and $\pi'$ will give rise to different contributions, because they now depend on $x$. 

This leads us to define 
\[ \phi_n[h](x):=\mathbb{E}\langle x|(M_h)^n|x\rangle .
\]
Note that $\phi_n[h]\!\!=\!\!\int\! \phi_n[h](x)dx$. Associating the label $x$ to the root of the tree $T_\bullet$ that corresponds to a partition $\pi$ we now have
\[z^{-n}\phi_n[h](x)=\sum_{T_\bullet \text{ with $n$ edges}} W(T_\bullet^x) .\] 
The weight $W(T_\bullet^x)$ of a tree $T_\bullet$ with root label $x$ is defined as follows~: Assign an integration variables $x_i$ to each black vertex, and assign $x$ to the black vertex that constitutes the root. Then assign the value $z^{-k}h(x_1)\cdots h(x_k)g_k(x_1,\cdots,x_k)$ to each white vertex whose neighbouring black vertices carry the variables $x_1,\cdots,x_k$. Finally, take the product over all vertices and integrates over all $x_i$ (except for the root $x$). By definition we set the tree consisting of a root without legs to one. Graphically the rules for the weights $W(T_\bullet ^x)$ are
\begin{center}
	\includegraphics[width=0.5 \linewidth]{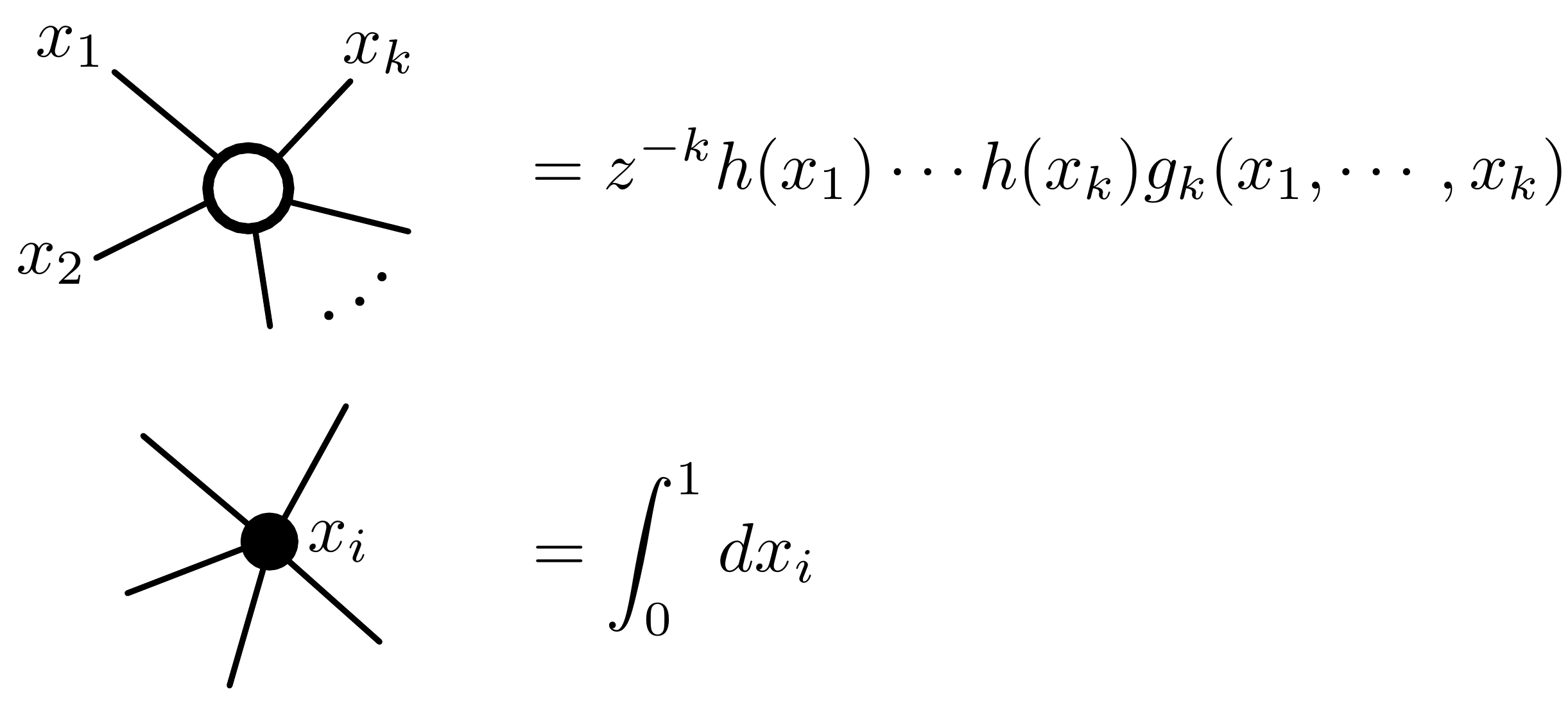}
\end{center}

Doing the sum over all $n$ is now easy~: We just relax the condition on the sum over trees with $n$ edges to trees of arbitrary size. We consider a generating function involving a sum over $\phi_n[h](x)$, 
\[
a_z(x):=\mathbb{E}\langle x|\frac{h}{z-M_h}|x\rangle=\frac{h(x)}{z}\sum_{n\ge 0} \frac{\phi_n[h](x)}{z^n}\overset{!}{=}\frac{h(x)}{z}\sum_{T_\bullet} W(T_\bullet^x),
\]
where the last equality is due to the correspondence with trees.

In order to establish the relation \eqref{eq:a-b} satisfied by $a_z(x)$ we consider the subset of trees $T_\circ$ whose root (still a black vertex) has a single leg only. This defines
\begin{equation}\label{eq:def_b}
	b_z(x):=\frac{z}{h(x)}\sum_{T_\circ} W(T_\circ^x) .
\end{equation}
Note that the weight $W(T_\bullet^x)$ of a tree whose root has $l$ legs is equal to the product of weights $W(T_{\circ,1}^x)\cdots W(T_{\circ,l}^x)$ of trees with a single leg on their root that arise by cutting the $l$ legs of $T_\bullet^x$. This implies 
\[\sum_{T_\bullet}W(T_\bullet^x)=\Big(1-\sum_{T_\circ} W(T_\circ^x)\Big)^{-1} ,\]
which yields the first relation in Eq.\eqref{eq:a-b}. 

For the second relation, we start with $T_\circ^x$ and cut the $l$ outgoing legs of the first white vertex. This generates a product of $l$ trees $T_{\bullet,i}^{x_i}$ whose weights satisfy 
\[W(T_\circ^x)=\frac{h(x)}{z}\int dx_1\cdots dx_l\, g_{l+1}(x,x_1,\cdots,x_l)W(T_{\bullet,1}^{x_1})\cdots W(T_{\bullet,l}^{x_l}).\]
Therefore, taking the sum over all trees $T_{\circ}$,
\begin{equation}
	b_z(x)=\sum_{l\ge0}\!\int \!\left(\prod_{i=1}^ldx_i  \frac{h(x_i)}{z}\!\sum_{T_\bullet}W(T_\bullet^{x_i})\right)g_{l+1}(x,x_1,\cdots,x_l).
\end{equation}
One recognizes the definition of $a_z(x_i)$ in this expression, which then implies the second relation in Eq.\eqref{eq:a-b}.

Both relations in Eq.~\eqref{eq:a-b} are the extremization conditions of the variational principle \eqref{eq:action}. As a last step we should therefore verify that $F[h]$ as defined in Eq.~\eqref{eq:F_expansion} coincides with the solution of the extremization problem from Eq.~\eqref{eq:action}. Here we show that their first derivates with respect to $h$ coincide for any $h$, as well as their value at $h=0$. 

Since $h(x)\,\delta \phi_n[h]/\delta h(x)=n\, \phi_n[h](x)$, one calculates from  Eq.\eqref{eq:F_expansion} that
\begin{equation}
	-h(x)\frac{\delta F[h](z)}{\delta h(x)}=\sum_{n\ge 1}\frac{\phi_n[h](x)}{z^{-n}}=\sum_{T_\bullet} W(T_\bullet^x)-1.
\end{equation}
Furthermore, one easily sees from our discussion of the multiplication of weights below Eq.\eqref{eq:def_b} that $a_z(x)b_z(x)=\sum_{T_\bullet} W(T_\bullet^x)-1$. This leads to 
\[-h(x)\frac{\delta F[h](z)}{\delta h(x)}=a_z(x)b_z(x)\]. 

On the other hand, starting from Eq.\eqref{eq:action}, one has\[ h(x)\frac{\delta F[h](z)}{\delta h(x)}=-\frac{h(x)b_z(x)}{z-h(x)b_z(x)}=-a_z(x)b_z(x),\] where we used Eq.~\eqref{eq:a-b} in the last line. Since $F[h=0](z)=\log(z)$ for both definitions \eqref{eq:action} and \eqref{eq:a-b}, the two expressions for $F[h](z)$ coincide.

\subsubsection{Proof using Kreweras duality}\label{subsec:proof_kreweras}
Let  $\tilde a_z$ be the resolvent with a marked variable $x$, such that $\int \tilde a_z(x)\,dx=G[h](z)$ is the resolvent, that is:
\[\tilde a_z(x) :=\mathbb{E}\langle x| \frac{ 1}{z- M_h} |x\rangle =\sum_{n\geq 0} z^{-n-1}\, \mathbb{E} \langle x| (M_h)^n|x\rangle  \]
From Eq.\eqref{eq:moment_free_cumulants}, the moments with marked variable $\phi_n[h](x):=\mathbb{E}\langle x|(M_h)^n|x\rangle$ are sum over non-crossing partitions, but without the integration over $x$.

Let $p_x$ be the part of $\pi$ containing $x$, and $p^*_x$ be the part of $\pi^*$ (the Kreweras dual) containing $x$ (we choose one of the two labelling of the edges of $\pi^*$ be naming them with their left (right) point when representing the partition as a loop connecting the points that belong to the same part). There are two ways to organise the sum over the number of points/edges and over their non-crossing partitions~: either by the cardinality of $p_x$ or that of $p_x^*$.

Let us first organise the sum by the cardinality of $p_x$. If $|p_x|=1$, then $x$ is a singlet in $\pi$. That is : we consider all marked partitions (with any number of points greater than one) such that $x$, the marked point, is a singlet. Summing over such partitions defines a function that we denote $\tilde b_z(x)$. That is,
\[ \tilde b_z(x):=\mathbb{E}\langle x| \frac{ z\,M_h}{z- M_h} |x\rangle^{[\mathrm{no}\, x]}=\sum_{n\geq 1} z^{1-n}\, \mathbb{E}\langle x| M_h^n|x\rangle^{[\mathrm{no}\, x]} 
, \]
where the "no $x$" upper script means that $x$ is not used in any of the intermediate indices in the product of matrices (i.e. when inserting a resolution of the identity).
If $k:=|p_x|\geq 2$, then the contribution of this partition to the product $\mathbb{E} \langle x| M_h^n|x\rangle$ splits into the product of $k$ contributions of the form $\mathbb{E}\langle x| M_h^{n_j}|x\rangle^{[\mathrm{no}\, x]}$, with $\sum_jn_j=n$. (Here we implicitly use the tree structure underlying the lattice of NC partitions). Since each NC partition appears only once, we get that
\[ \tilde a_z(x) = \sum_{k\geq 0} z^{-1-k}\, [\tilde b_z(x)]^k = \frac{1}{z- \tilde b_z(x)} .\]
This is the first relation in \eqref{eq:a-b} if we define $a_z=h\,\tilde{a_z}$ and $b_z=\tilde b_z/h$.

Let us now organise the sum by the cardinal of $p_x^*$. Actually we shall organise the sum involved in $\tilde b_z(x)$ (with no repetition of $x$). Let $k:=|p_x^*|\geq 1$. The contributions of the marked NC partitions with $|p_x^*|=k$ will each involve a factor $g_{k+1}(x,x_1,\cdots,x_k)$. Using again the tree structure underlying the lattice of NC partitions, we then read that\footnote{In the discrete level, because of the "$\mathrm{no}\, x$" constraint in the definition of $b(x)$, the indices to sum over, representing the integrals over $x_1,\cdots, x_k$, should be different from the marked point $x$. Similarly, the functions $  \tilde a_z(x_j)$ which arise in the relation $ \tilde b_z(x) = R_{0}[\tilde a_z](x)$ should actually be $\tilde a_z(x_j)^{\mathrm{no}\,x}$, involving the matrix element $\langle x_j|M_h^n|x\rangle^{\mathrm{no}\,x_j}$. But this does not matter in the continuum limit, since all the functions to be integrated over are smooth and $x\not= x_j$ in the integration.}
\[ \tilde b_z(x) = \sum_{k\geq 0} h(x) \int dx_1\cdots dx_k\, g_{k+1}(x,x_1,\cdots,x_k) h(x_1) \tilde a_z(x_1)\cdots h(x_k) \tilde a_z(x_k) .\]
Using the definition for $ R_{0}[p](x)$ below Eq.\eqref{eq:a-b} this reads
\[  \tilde b_z(x) = h(x)\,  R_{0}[h\, \tilde a_z](x).\]
With $a_z=h\,\tilde{a_z}$ and $b_z=\tilde b_z/h$, this is the second relation in Eq.\eqref{eq:a-b}

Finally note that,
\[ h(x)\frac{\delta F[h](z)}{\delta h(x)} = 1 -z  \tilde a_z(x) = - \tilde a_z(x) \tilde b_z(x)=a_z(x)b_z(x).\]
Together with the boundary value, $F[0](z)=\log z$, this fixes $F[h](z)$ (at least as formal series in $h$).

\subsubsection{Proof using operator valued free probability}\label{subsec:proof_op-val}
This section recalls some basic definitions of operator-valued free probability theory and shows how the relation between the R- and the Cauchy-transform (Theorem \ref{thrm:R-Cauchy}) can be used to deduce our main result (Theorem \ref{thrm:F}). Of course, the relation between $R$- and Cauchy transform also uses implicitly the tree structure of non-crossing partitions. We closely follow \cite[chpt. 9]{Mingo2017Free} and \cite{Shlyakhtenko1998Gaussian} and start with the definition of the operator-valued moments for a general unital algebra $\mathcal A$ which later becomes the matrix algebra formed by the matrices $M$.
\vspace{5pt}
\begin{definition}
	Let $\mathcal A$ be a unital algebra and consider a unital subalgebra $\mathcal D \subset \mathcal A$. Then $E^\mathcal{D}:\mathcal A \to \mathcal D$ is called a \emph{conditional expectation value} (with amalgamation over $\mathcal D$) if for all $a\in \mathcal A$ and $d, d'\in \mathcal D$ one has $E^\mathcal{D}[d]\in\mathcal D$ and $E^\mathcal{D}[d a d']=dE^\mathcal{D}[a]d'$. 
	
	For any choice of $d_1,\cdots,d_{n-1}\in\mathcal D$, the \emph{operator-valued} (or \emph{$\mathcal{D}$-valued}) \emph{moments} of $a$ are defined as $E^\mathcal{D}[a d_1a\cdots ad_{n-1}a]\in\mathcal D$ and the collection of all operator-valued moments define the operator-valued distribution of $a$.
\end{definition}

We will now consider the special case where the elements $M\equiv a\in\mathcal A$ are random matrices of size $N$ satisfying properties (i)-(iii), and the elements $\Delta\equiv d\in\mathcal D$ are diagonal matrices of size $N$. Note that $\mathcal{D}$ is indeed a subalgebra of $\mathcal A$ and that in the large $N$ limit we have $\mathcal D\to L^\infty[0,1]$. We also define explicitly a conditional expectation value adapted to our choice of $\mathcal D\subset\mathcal A$. For $M\in \mathcal{A}$,
\begin{equation}\label{eq:cond_expectation}
	E^\mathcal{D}[M]:=\mathrm{diag}(\mathbb E[M_{11}],\cdots,\mathbb E[M_{NN}]).
\end{equation}
Clearly, this definition satisfies the defining properties of a conditional expectation value.

As in scalar free probability, one can define operator-valued free cumulants as follows.
\begin{definition}
	The \emph{$\mathcal D$-valued free cumulants} $\kappa^\mathcal{D}_n:\mathcal A^n\to\mathcal{D}$ are implicitly defined by 
	\begin{equation}\label{eq:op-val-free-cumulant}
		E^\mathcal{D}[M_1\cdots M_n]=:\sum_{\pi\in NC(n)} \kappa_\pi^\mathcal{D}(M_1,\cdots,M_n)
	\end{equation}
	where $\kappa_\pi^\mathcal{D}$ is obtained from the family of linear functions $\kappa_n^\mathcal{D}:=\kappa_{1_n}^\mathcal{D}$ by respecting the nested structure of the parts appearing in $\pi$ as explained in the following example.
\end{definition}

\begin{example}\label{ex:nested}
	For $\pi=\{\{1,3\},\{2\},\{4,5\},\{6\}\}$, which corresponds to the dotted lines in the following figure, $\kappa_\pi^\mathcal{D}$ is defined as
	\begin{center}
		\includegraphics[width=0.2 \linewidth]{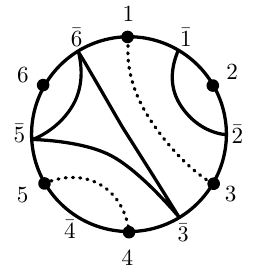}
	\end{center}
	\begin{equation*}
		\kappa_\pi^\mathcal{D}(M_1,M_2,M_3,M_4,M_5,M_6):=\kappa_2^\mathcal{D}(M_1\cdot \kappa_1^\mathcal{D}(M_2),M_3)\cdot \kappa^\mathcal{D}_2(M_4,M_5)\cdot \kappa^\mathcal{D}_1(M_6).
	\end{equation*}
	Note that one deals with matrix products, specifically emphasized by the dot $\cdot$ in this example, which is omitted elsewhere.
\end{example}
\medskip
Next we would like to relate $\kappa_n^\mathcal{D}$ to the local free cumulants $g_n$. In the large $N$ limit with $x=i/N$, we introduce the notation $E^\mathcal{D}[M](x):=E^\mathcal{D}[M]_{ii}\in \mathbb R$ to denote a diagonal elements of $\mathcal D$. By Eq.~\eqref{eq:moment_free_cumulants}, we can express the $\mathcal D$-valued moments as
\begin{equation*}
	E^\mathcal{D}[M\Delta_1\cdots M\Delta_n M](x)=\!\!\!\!\! \sum_{\pi\in NC(n+1)}\int \mathrm d \vec x^{(n)}\Delta_1(x_1)\cdots \Delta_n(x_n)g_\pi(\vec x^{(n)},x) \delta_{\pi^*}(\vec x^{(n)},x)
\end{equation*}
Here we interchanged the roles of $\pi$ and $\pi^*$ which does not change the sum. Comparing to the definition of operator-valued free cumulants, this suggest the following identification.
\begin{proposition} Let $\pi\in NC(n+1)$, we have
	\begin{equation}\label{eq:op-val-cum_pi}
		\kappa_\pi^\mathcal{D}(M\Delta_1,\cdots,M\Delta_n,M)(x)=\int \mathrm d \vec x^{(n)} \Delta_1(x_1)\cdots \Delta_n(x_n) g_\pi(\vec x^{(n)},x) \delta_{\pi^*}(\vec x^{(n)},x)
	\end{equation}
\end{proposition}
\begin{proof}
	We must check that this identification can be consistently obtained from the case $\pi=1_{n+1}$ by respecting the nested structure appearing in $\kappa_\pi^\mathcal{D}$. That is, we define
	\begin{equation}\label{eq:op-val-cum_n}
		\kappa_{n+1}^\mathcal{D}(M\Delta_1,\cdots,M\Delta_n,M)(x):=\int \mathrm d \vec x^{(n)} \Delta_1(x_1)\cdots \Delta_n(x_n) g_{n+1}(\vec x^{(n)},x),
	\end{equation}
	and show that this implies the proposition. It is important to have included the diagonal $\Delta_i$'s in this definition, since this allows us to resolve nested terms such as $\kappa_2^\mathcal{D}(M\, \kappa_1^\mathcal{D}(M),M)$. In fact, one soon notices that Eq.\eqref{eq:op-val-cum_pi} is precisely the definition of the nested structure of $\kappa_\pi^\mathcal{D}$. 
	
	We illustrate this using the above example with $\pi=\{\{1,3\},\{2\},\{4,5\},\{6\}\}$ and Kreweras complement $\pi^*=\{\{\bar 1,\bar 2\},\{\bar 3,\bar 5,\bar 6\},\{\bar 4\}\}$.
	The definition of $\kappa_{n+1}^\mathcal{D}$ implies that the l.h.s of Eq.\eqref{eq:op-val-cum_pi} becomes
	\begin{align*}
		& \kappa_2^\mathcal{D}(M\Delta_1 \kappa_1^\mathcal{D}(M\Delta_2),M\Delta_3)(x)\, \kappa^\mathcal{D}_2(M\Delta_4,M\Delta_5)(x)\, \kappa^\mathcal{D}_1(M)(x)\\
		& = \int \mathrm d x_1\,\Delta_1(x_1)\Delta_2(x_1)g_1(x_1) \Delta_3(x)g_2(x_1,x) \int \mathrm d x_4\, \Delta_4(x_4)\Delta_5(x)g_2(x_4,x) g_1(x)
	\end{align*}
	This corresponds indeed to the r.h.s. where $\delta_{\pi^*}(x_1,\cdots,x_6)=\delta(x_1-x_2) \delta(x_3-x_5)\delta(x_5-x)$ and we identified $x_6\equiv x$. An arbitrary $\pi\in NC(n+1)$ can be tackled in the same way identifying $x_{n+1}\equiv x$.
\end{proof}

This result also explains how the structure of $g_\pi \delta_{\pi^*}$ which we encountered in Eq.\eqref{eq:moment_free_cumulants} fits into the free probability picture. Earlier, we could only ascertain that the family $\tilde \kappa_\pi:=\int g_{\pi^*}(\vec x)\,\delta_\pi(\vec x)d\vec x$ are not the (scalar) free cumulants of $M$ because they where not multiplicative ($\tilde\kappa_\pi\tilde \kappa_\sigma\neq\tilde \kappa_{\pi\cup\sigma}$). Now we understand that they are nonetheless free cumulants, but in the operator valued setting with amalgamation over diagonal matrices. More precisely $\tilde \kappa_\pi =\ntr\left( \kappa_\pi^\mathcal{D}(M,\cdots,M)\right)$. This also suggests that calling the family of functions $g_n$ "local free cumulants" seems to be a good name choice.
\smallskip
\begin{definition} \label{def:op-val-R-and-Cauchy}
	The \emph{$\mathcal{D}$-valued $R$-transform}, $R_M:\mathcal D\to\mathcal D$ of an element $M\in\mathcal A$ is defined by
	\begin{equation}
		R_M(\Delta) :=\sum_{n\ge0} \kappa_{n+1}^\mathcal{D}(M\Delta,\cdots,M\Delta,M)
	\end{equation}
	and the \emph{$\mathcal{D}$-valued Cauchy transform} (or \emph{resolvent}) $G_M:\mathcal{D}\to\mathcal D$ is defined by
	\begin{equation}
		G_M(\Delta)=E^\mathcal{D}[\frac 1 {\Delta-M}]=\sum_{n \ge 0} E^\mathcal{D}[\Delta^{-1}(M\Delta^{-1})^n]
	\end{equation}
\end{definition}
\begin{theorem}[see Thrm. 11 in Chpt. 9 of \cite{Mingo2017Free}]\label{thrm:R-Cauchy}	
	Similarly to the scalar-valued case, here the $R$- and Cauchy transforms are related by
	\begin{equation}
		G_M(\Delta)=\frac{1}{\Delta-R_M(G_M(\Delta))}.
	\end{equation}
\end{theorem}
The $\mathcal{D}$-valued Cauchy transform can be related to its scalar analogue $G(z)$ (denoted without the subscript) by 
\begin{equation*}
	G(z):=\ntr(G_M(z\mathbb I))=\E\big[\ntr\big(\frac{1}{z\mathbb I-M}\big)\big].
\end{equation*}
Let us now consider $M_h=h^{1/2}Mh^{1/2}\in\mathcal A$ and define (with $x=i/N$)
\begin{equation}
	\tilde a_z(x):=\lim_{N\to\infty}G_{M_h}(z\mathbb I)_{ii}.
\end{equation}
The scalar Cauchy transform of $M_h$ is then $G(z)=\int_0^1 \mathrm d x \, \tilde a(x)$.
Furthermore, from Eq.\eqref{eq:op-val-cum_n} one sees that
\begin{equation}\label{eq:op-val-R}
	R_{M_h}(\Delta)(x)=\sum_{n\ge 0} \int \mathrm d\vec x^{(n)} \, \Delta(x_1)h(x_1)\cdots \Delta(x_n)h(x_n)h(x)g_{n+1}(\vec x^{(n)},x)
\end{equation}
Together with Theorem \ref{thrm:R-Cauchy} we therefore obtain,
\begin{equation}
	\tilde{a}_z(x)=\frac 1 {z-R_{M_h}(\tilde a_z)(x)}=\frac{1}{z-h(x)R_0[h\tilde a_z](x)}.
\end{equation}
In the last equality we used the definition of $R_0$ from Eq.\eqref{eq:R_0}. Redefining $a_z(x)=h(x)\tilde a_z(x)$ we obtain the extremization conditions in Eq.\eqref{eq:a-b}, which are equivalent to the variational principle in Theorem \ref{thrm:F}.

\subsection{No free multiplicative convolution}\label{subsec:no_free_mulit_convolution}

In the case of matrices rotated by Haar random unitaries we will see in subsection \ref{subsec:Haar} that the spectral measure $\sigma_I$ of $M_h$ is related to that of $h$ and $M$, respectively denoted by $\nu$ and  $\sigma$, via a free multiplicative convolution as $\sigma_I=\nu \boxtimes  \sigma$. This is, in fact, always true if $M$ and $h$ are free in the sense of free probability. For $M$ a Haar-randomly rotated matrix, this is the case~: $h$ and $M$ are free. Here we will show that we cannot obtain $\sigma_I$ by free multiplicative convolution of $\nu$ and $\sigma$ in the general case of structured random matrices where the local free cumulants $g_n$ are not constant. In return, this means, that structured random matrices $M$ are not free from deterministic diagonal matrices $h$.

Let $S_I$ (resp. $S_0$) be the $S$-transform of the measure $\sigma_I$ (resp. $\sigma$). Recall that $S_I(w)=\frac{w+1}{wz_I}$ with $w+1= z_IG_I(z_I)$ and similarly  $S_0(w)=\frac{w+1}{wz_0}$ with $w+1= z_0G_0(z_0)$ (see appendix \ref{app:glossary} for a definition of the S-transform). Hence, if $\sigma_I=\nu \boxtimes  \sigma$ so that $S_I(w)=S_0(w)S_h(w)$, we should have
\[ \frac{z_0}{z_I}= S_h(w) .\]
In particular, if $\sigma_I=\nu \boxtimes  \sigma$ then ${z_0}/{z_I}$ is independent of the distribution of $M$, that is, it is independent of the local free cumulants $g_n$.
We check below that ${z_0}/{z_I}$ actually depends on the local free cumulants $g_n$, so that $\sigma_I\not=\nu \boxtimes  \sigma$.

The equation for $z_I$ reads
\[ w = \int\!dx\frac{h(x)b_z(x)}{z_I-h(x)b_z(x)} =  \int\!dx\left(\frac{1}{z_I} h(x)b_z(x) + \frac{1}{z_I^2} (h(x)b_z(x))^2+\cdots \right).\]
This is an equation for $1/z_I$. We have to take into account that $b_z(x)$ also depends on $z$. Let us write $b_z(x)=b_1(x) +\frac{1}{z_I} b_2(x) +\cdots$. We have $b_1(x)=g_1(x)$ and $b_2(x)=\int\! dy\, g_2(x,y)h(y)$. Solving for $1/z_I$ we get
\[
\frac{1}{z_I} = \frac{w}{[hb_1]}\left( 1- \frac{ [(hb_1)^2] + [hb_2]}{[hb_1]^2}w+\cdots \right)
\]
where $[\cdots]$ is short notation for integration over $x$, i.e. $[f]=\int\!dxf(x)$. The formula for $z_0$ is obtained from that of $z_I$ be setting $h=1$. Thus
\[
\frac{z_0}{z_I}= \frac{[g_1]}{[hg_1]}(1+O(w)),
\]
and ${z_0}/{z_I}$ depends on the local free cumulants $g_n$, and $S_I(w)\not= S_h(w)S_0(w)$ or equivalently $\sigma_I\not=\nu \boxtimes  \sigma$.

Remark that, to next order in $w$ we have
\[
\frac{z_0}{z_I}= \frac{[g_1]}{[hg_1]}\left(1- (\frac{ [(hg_1)^2]}{[hg_1]^2}-\frac{ [(g_1)^2]}{[g_1]^2})w  - (\frac{ [[g_2(h\times h)]]}{[hg_1]^2}-\frac{ [[g_2]]}{[g_1]^2})w +O(w^2)\right)
\]
For constant local free cumulants (i.e. for unstructured matrices), $g_1\leadsto \kappa_1$ and $g_2\leadsto \kappa_2$, the r.h.s. of the previous equation is independent of those cumulants (it depends only on $h$). More precisely the second term proportional to $w$ vanishes since $[[g_2(h\times h)]]\leadsto \kappa_2 [h]^2$ and $[hg_1]\leadsto \kappa_1[h]$, while the remaining terms become
\[
\frac{z_0}{z_I}= \frac{1}{[h]}\left(1- \frac{ [h^2]-[h]^2}{[h]^2}w+O(w^2)\right)
\]
This coincides with the expansion of $S_h(w)$ as expected.

\section{Applications} \label{sec3}
In this section we apply the formulae (\ref{eq:resolvent}-\ref{eq:a-b}) to some explicit random matrix ensembles. For Wigner matrices we show that this produces the well known Wigner semi-circle law. We also show how to generalize this to the structured case where the variance of diagonal entries can vary. For matrices rotated by Haar random unitaries, we show that our method reduces to free multiplicative convolution when interested in the spectrum of subblocks. Finally we apply our method to the stationary distribution of the Quantum Symmetric Simple Exclusion Process (QSSEP), a structured random matrix ensembles for which the local free cumulants are known. This application to QSSEP is new and extends the results previously presented in \cite{Bernard2023Exact}.

\subsection{Wigner matrices} \label{subsec:wigner}
Wigner matrices are characterized by the vanishing of its associated free cumulants of order strictly bigger than two. Thus, for Wigner matrices only $g_1$ and $g_2$ are non vanishing and both are $x$-independent. All $g_n$, $n\geq 3$, are zero. Without loss of generality we can choose $g_1=0$ and we set $g_2=s^2$. Then $F_0[p]=\frac{s^2}{2}\int \!dx dy\, p(x)p(y)$ and $R_0[p]=s^2 \int\! dx\, p(x)$. For the whole interval $h(x)=1$ (considering a subset will be equivalent), the saddle point equations become
\[
a = \frac{1}{z- b},\ b = s^2\, A, 
\]
with $A=\int \!dx\, a(x)$. This yields a second order equation for $A$, i.e. $A^{-1}=z-s^2 A$. Solving it, with the boundary condition $A\sim \frac{1}{z}+\cdots$ at $z$ large, gives
\[ 
A = \frac{1}{2s^2}\left( z - \sqrt{z^2-4s^2}\right)
\]
Thus the cut is on the interval $[-2s,+2s]$ and the spectral density is
\begin{equation}
d\sigma(\lambda) = \frac{d\lambda}{2\pi s^2}\, \sqrt{4s^2-\lambda^2}\ 1_{\lambda\in[-2s,+2s]}
\end{equation}
Of course, that's Wigner's semi-circle law.

%\paragraph{Inhomogeneous Variance.}
%Next we consider $N\times N$ Wigner matrices with zero mean and variance
%\[
%\mathbb{E}[M_{ij}M_{kl}] = \frac{1}{N} \delta_{jk}\delta_{il}\, g_2(\frac{i}{N},\frac{j}{N}) ,
%\]
%with $g_2(x,y)$ a (smooth) function. It is clear that the three fundamental properties (i)-(iii) are satisfied. We restrict to diagonal covariances $g_2(x,y)=s^2(x)\delta(x-y)$, because otherwise we cannot find closed expressions for the spectrum. The saddle point equation is then a quadratic equation for $a_z(x)$ which, in the case $h(x)=1$, reads $a_z(x)(z-s^2(x)a_z(x))=1$ so that
%\[
%a_z(x) = \frac{1}{2s^2(x)}(z-\sqrt{z^2-4s^2(x)}) .
%\]
%The resolvent is $G(z)=\int\! dx\, a_z(x)$. Its discontinuity at the cut is the sum of the discontinuities for each value of $x$. This yields for the spectral density
%\begin{equation}
%d\sigma(\lambda) = \frac{d\lambda}{2\pi}\int\! \frac{dx}{s^2(x)}\, \sqrt{4s^2(x)-\lambda^2}.
%\end{equation}

\subsection{Haar-randomly rotated matrices} \label{subsec:Haar}

We consider matrices of the form $M=UDU^\dag$, with $U$ Haar distributed over the unitary group and $D$ a diagonal matrix with spectral density $\sigma$ in the large $N$ limit. For such matrices, it is known that the local free cumulants are constant and equal to the free cumulants of $\sigma$, that is
\begin{align} \label{eq:Haar-loop}
	g_n(\vec x):=\lim_{N\to\infty}N^{n-1}\mathbb{E}[M_{i_{1}i_{2}}M_{i_{2}i_{3}}\cdots M_{i_{n}i_{1}}]^c=\kappa_n(\sigma) .
\end{align}
The proof resorts to the HCIZ integral and is outlined in appendix \ref{app:free_cum_haar}.

\paragraph{Spectrum of $M$ (consistency check).}
Of course the spectrum of the whole matrix $M$ is that of $D$ with spectral density $\sigma$. Let us check this within our approach via Theorem \ref{thrm:F}. With  $g_n=\kappa_n(\sigma)$ and $h(x)=1$ (we consider the whole matrix $M$), Eqs.\eqref{eq:a-b} become 
\[
A =\frac{1}{z- b_z(A)},\  b_z(A)= \sum_{k\geq 1} A^{k-1}\kappa_k(\sigma), %:=R_{\sigma}(A_z)
\]
with $A=\int\!dx\, a_z(x)$. Let us recall some basics definition from free probability. For any measure $\sigma$ of some random variable $X$, let $G_\sigma(z)=\mathbb{E}[\frac{1}{z-X}]= \sum_{n\geq 0}z^{-n-1}m_n(\sigma)$  and $K_\sigma(z)=\sum_{n\geq 0}z^{n-1}\kappa_n(\sigma)$, with $m_n$ and $\kappa_n$ the $n$-th moments and free cumulants, respectively. As well known from free probability, $G_\sigma$ and $K_\sigma$ are inverse functions, i.e. $K_\sigma(G_\sigma(z))=z$. Comparing with the previous equation, we see that $b_z(A)=K_{\sigma}(A)-A^{-1}$. The equation $A=1/(z- b_z(A))$ can thus be written as $z= b_z(A)+A^{-1}=K_{\sigma}(A)$, and hence 
\[
A=G_{\sigma}(z)
\]
As a consequence, the resolvent of $M$ is equal to $G_{\sigma}(z)$ and the spectral density of $M$ is indeed that of $D$, as it should be.

\paragraph{Spectrum of a subblock reduces to free compression.}
We now consider an interval $I=[0,\ell]\in[0,1]$ and compute the spectrum $d\sigma_I(\lambda)$ of the subblock of $M$ with $\ell N$ rows and columns that corresponds to this interval, i.e. $h(x)=1_{x\in I}$. The saddle point equations \eqref{eq:a-b} impose $b_z(x)$ to be independent of $x$ and $a_z(x)=0$ for $x\not\in I$. They then read 
\begin{equation*}
	b_z = \sum_{k\geq 1} A_\ell^{k-1}\kappa_k(\sigma),\quad A_\ell=\frac{\ell}{z-b_z}
\end{equation*}
with $A_\ell=\int_{I}\!dx\, a_z(x)$. These two equations imply
$z= K_{\sigma}(A_\ell) -( 1-\ell)/{A_\ell}$.

Let us now define (following a remark by Ph. Biane) the freely-compressed measure $\sigma^{(t)}$ defined from $\sigma$ by compressing its free cumulants by a factor $1/t$, that is 
\begin{equation*} %\label{eq:free_dilation}
	\kappa_k(\sigma^{(t)}) :={t}^{-1} \kappa_k(\sigma).
\end{equation*}
We have $K_{\sigma^{(t)}}(w) = \frac{1}{t} K_{\sigma}(w) + (\frac{t-1}{t}) \frac{1}{w}$.
The equation for $A_\ell$ thus reads $K_{\sigma^{(t)}}(A_\ell)=\frac{z}{\ell}$ and hence
\begin{equation*}
	A_\ell(z) = G_{\sigma^{(\ell)}}(\frac{z}{\ell}) .
\end{equation*}
This implies, $\int \frac{d\sigma_I(\lambda)}{z-\lambda} = \int \frac{d\sigma^{(\ell)}(X)}{z-\ell X}$, so that
\begin{equation}
	d\sigma_I(\ell \lambda) = d\sigma^{(\ell)}({\lambda})
\end{equation}
That is~: the spectral measure of a subblock $I=[0,\ell]$ of $M$ of relative size $\ell$ is that of the freely-compressed measure $\sigma^{(\ell)}$ but for the compressed eigenvalue ${\lambda}/{\ell}$.

\paragraph{Spectrum via free multiplicative convolution.}
For illustrating purposes, we present an explicit derivation of the well-known result that for any initial spectral measure $\sigma$ -- that of the diagonal matrix $D$ -- and any function $h$ -- defining the total spectral measure $\nu$ --, the spectral measure of $h^{1/2}Mh^{1/2}$ is the free multiplicative convolution of $\nu$ and $\sigma$, that is
\begin{equation}
	\sigma_I^\mathrm{tot} = \nu \boxtimes \sigma.
\end{equation}
Note that $\sigma_I^\mathrm{tot}$ includes potential zero eigenvalues in the case where $h(x)=1_{x\in I}$ is the indicator function on an interval. In contrast to this, $\sigma_I$, which we obtained in the last subsection via free compression as the spectrum of a subblock $M_I\subset M$, does not contain these zero eigenvalues. The two spectra are related by Eq.\eqref{eq:total_spectrum}.

Let $S_0$, $S_h$, $S_I$ (resp. $G_0$, $G_h$, $G_I$) the S-transform (resp. the resolvent) of the spectral measure $\sigma$, $\nu$, $\sigma_I^\mathrm{tot}$.
Recall the relation between the $S$-transform and the $R$-transform as $C(zS(z))=z$ with $C(z)=zR(z)$. Recall also the relation between the resolvent and the $S$-transform: $S(zG(z)-1)=G(z)/(zG(z)-1)$. Setting $w=zG(z)-1$, it can be written as $S(w)=\frac{w+1}{zw}$ with $z$ implicitly depending on $w$ via $w+1=zG(z)$.

We have to prove $S_I(w)=S_h(w)S_0(w)$.

Let $w:=zG_I(z)-1$, so that $S_I(w)=\frac{w+1}{zw}$ with $zG_I(z)=w+1$. 
The resolvent $G_I$ is defined by the saddle point equations
\begin{equation*}
	G_I(z) = \int\!\frac{dx}{z-h(x) R_0(A)},\quad \mathrm{with}\ A= \int\! \frac{h(x)\,dx}{z-h(x)R_0(A)} .
\end{equation*}
Using these relations we have $w= R_0(A)\int\! \frac{h(x)\,dx}{z-h(x)R_0(A)}=AR_0(A)$ thus $w=C_0(A)$, and hence $A=wS_0(w)$. As a consequence, the relation $A= \int\! \frac{h(x)\,dx}{z-h(x)R_0(A)}$ can be written as (using $A=wS_0(w)$)
\begin{equation*}
	w= \int\! \frac{h(x)dx}{zS_0(w)-h(x)}
\end{equation*}
Let $u=zS_0(w)$, then (using $S_I(w)=\frac{w+1}{zw}$)
\begin{equation*}
	S_I(w)= S_0(w)\, \frac{w+1}{uw},\quad \mathrm{with}\ w+1 = \int\!\frac{u\,dx}{u-h(x)}
\end{equation*}
Now $\frac{w+1}{uw}=S_h(w)$, because for $h$ diagonal, $G_h(u)=\int\!\frac{dx}{u-h(x)}$ and hence $S_h(w)=\frac{w+1}{uw}$ with $w+1=\int\!\frac{u\,dx}{u-h(x)}$.
Thus $S_I(w)=S_0(w)S_h(w)$ or equivalently $\sigma_I^\mathrm{tot} = \nu \boxtimes \sigma$.

\subsection{QSSEP}\label{subsec:QSSEP}
The open Quantum Symmetric Simple Exclusion Process (QSSEP) is a quantum stochastic process that is supposed to model diffusive transport in one dimensional chaotic many-body quantum systems in the mesoscopic regime \cite{Bernard2019Open,Hruza2023Coherent}. Mathematically it is a one-dimensional chain with $N$ sites occupied by spinless free fermions $c^\dagger_j$ with noisy hopping rates and coupled to boundary reservoirs at $j=1,N$ that inject and extract fermions with rates $\alpha_{1,N}$ and $\beta_{1,N}$, respectively.
The key quantity of interest is the matrix of coherences with elements $M_{ij}(t):=\mathrm{Tr}(\rho_{t}\,c_{i}^{\dagger}c_{j})$ which contains all information about the system since the evolution of QSSEP preserves Gaussian fermionic states \cite{Bernard2019Open}. The $N\times N$ matrix $M$ undergoes a stochastic evolution of the form
\begin{equation} \label{eq:M_evolution}
	dM(t)=i[dh_t,M(t)]-\frac 1 2 [dh_t,[dh_t,M(t)]]+\mathcal L[M] dt
\end{equation}
with 
\begin{align*}
	dh_t=
	\begin{pmatrix}
		0 				&dW_t^1 	& 					& 			\\
		d\overline W_t^1	&\ddots 	&\ddots 			& 			\\
		&\ddots		&\ddots				& dW_t^{N-1}\\
		&			&d\overline W_t^{N-1}& 0			\\
	\end{pmatrix}
\end{align*}
where $dW_t^j:=W_{t+dt}^j-W_{t}^j$ are the increments of complex Brownian motions, independent for each site $j$, and 
\begin{equation*}
	\mathcal{L}[M]_{ij}=\sum_{p\in{1,N}}(\delta_{pi}\delta_{pj}\alpha_p - \frac 1 2 (\delta_{ip}+\delta_{jp})(\alpha_p+\beta_p)M_{ij})
\end{equation*}
The stochastic evolution has a unique stationary distribution \cite{Bernard2021Solution} that is characterized by its local free cumulants as \cite{Biane2021Combinatorics} (again with the notation $\vec x=(x_1,\cdots,x_n)$)
\begin{equation}\label{eq:g_QSSEP}
	\sum_{\pi\in NC(n)} g_\pi(\vec x)=\min(\vec x)=:\varphi_n(\vec x).
\end{equation}
The functions $g_n$ can be viewed as the free cumulants of the indicator functions $\mathbb{I}_{x}(y):=1_{y<x}$ with respect to the Lebesgue measure, since the moments of these functions are precisely $\mathbb{E}[\mathbb{I}_{x_{1}}\cdots\mathbb{I}_{x_{n}}]=\min(\vec x)$.

From the point of view of physics we are interested in the spectra of $M$ and of its subblocks in order to compute the entanglement entropy in the QSSEP \cite{Bernard2023Exact}. 

\paragraph{Explicit expression for $F_{0}[p]$.}

As a first step, we compute the function $F_0[p]$ that contains the initial data through the functions $g_n$. Defining $\mathbb{I}_{[p]}(y):=\int_{y}^{1}dx\,p(x)$, we shall prove that for QSSEP,
\begin{align} \label{eq:R0-implicit}
	F_{0}[p] =w-1-\int_{0}^{1}dx\log[w-\mathbb{I}_{[p]}(x)] &,\ \mathrm{with}\ 
	\int_{0}^{1}\frac{dx}{w-\mathbb{I}_{[p]}(x)}=1,
\end{align}

We define the free cumulant and the moment generating function,
\begin{align*}
	K_{[p]}(w) =\sum_{n\ge0}w^{n-1}g_{n}[p], \quad 
	G_{[p]}(w) =\sum_{n\ge0}w^{-n-1}\varphi_{n}[p],
\end{align*}
where $\varphi_{n}[p] :=\int d\vec x\, \varphi(\vec x)p(x_{1})\cdots p(x_{n})$ and $g_{n}[p] :=\int d\vec x\, g(\vec x)p(x_{1})\cdots p(x_{n})$. By convention, we set $g_{0}[p]=\varphi_{0}[p]\equiv1$. We have $\varphi_{n}[p]=\sum_{\pi\in NC(n)}g_{\pi}[p]$. By results from free probability theory, these two functions are inverses of each other, $K_{[p]}(G_{[p]}(w)) =w$. Integrating Eq.\eqref{eq:g_QSSEP} and using $\mathbb{I}_{[p]}(y)=\int_{0}^{1}dx\,p(x)\mathbb{I}_{x}(y)$, the Cauchy transform $G_{[p]}$ can be written as
\begin{equation}\label{eq:z_definition}
	G_{[p]}(w)=\int_{0}^{1}\frac{dx}{w-\mathbb{I}_{[p]}(x)}.
\end{equation}
Since the initial data function $F_0$ is such that $F_{0}[vp]=\sum_{n\ge1}\frac{v^{n}}{n}g_{n}[p]$, we have $1+v\partial_{v}F_{0}[vp]=vK_{[p]}(v$).
Define now a new variable $w$, depending on $v$ and $p$, such that $v=G_{[p]}(w)$. Using $K_{[p]}(G_{[p]}(w))=w$, the equation $1+v\partial_{v}F_{0}[vp]=vK_{[p]}(v)$ then becomes 
\[
1+v\partial_{v}F_{0}[vp]=vw
\]
Integrating w.r.t.\ $v$ yields (with the appropriate boundary condition $F_0[0]=0$)
\[
F_{0}[vp]=vw-1-\int_{0}^{1}dx\log[v(w-\mathbb{I}_{[p]}(x))]
\]
Indeed, computing the $v$-derivative of the l.h.s gives $v\partial_vF_{0}[vp]=vw-1+(\frac{\partial w}{\partial v})[v-\int_0^1\frac{dx}{w-\mathbb{I}_{[p]}(x)}]$ which, using equation \eqref{eq:z_definition}, becomes $v\partial_vF_{0}[vp]=vw-1$. Setting $v=1$ one obtains Eq.\eqref{eq:R0-implicit}.

\paragraph{Differential equation for $b_z(x)$.}
Next we derive a differential equation for $b_z(x)$, see Eq.\eqref{eq:b-diff} below.

Using Eq.\eqref{eq:R0-implicit}, the relation $b_z(x)=\frac{\delta F_0[a_z]}{\delta a_z(x)}$ becomes 
\begin{equation*}
	b_z(x)= \int_0^x \frac{dy}{w-\mathbb{I}_{[a_z]}(y)},\ \mathrm{with}\ \int_{0}^{1}\frac{dx}{w-\mathbb{I}_{[a_z]}(x)}=1.
\end{equation*}
Thus, $b_z(x)$ satisfies the boundary conditions $b_z(0)=0$ and $b_z(1)=1$. Furthermore, $1/b_z'(x)=w-\mathbb{I}_{[a_z]}(x)$ and $(1/b_z'(x))'=a_z(x)$.
Using $a_z(x)=\frac{h(x)}{z-h(x)b_z(x)}$ from the saddle point equation gives, after some algebraic manipulation,
\begin{equation}\label{eq:b-diff}
	z b''(x) + h(x)(b'(x)^2-b(x) b''(x))=0 .
\end{equation}

For $h(x)=1_{x\in I}$, that is $h(x)=0$ for $x\not\in I$ and $h(x)=1$ for $x\in I$, this yields
\begin{equation}\label{eq:b_z_equation}
	\begin{cases}
		[\log(z-b_z(x))]''=0, & \text{ if }x\in I \\
		b_z(x)''=0, & \text { if } x\notin I
	\end{cases}
\end{equation}
with boundary conditions $b_z(0)=0$ and $b_z(1)=1$.

\paragraph{Spectrum of $M$.}
First we present the derivation of the easier case where $I=[0,1]$. In this case a solution of Eq.\eqref{eq:b_z_equation} with correct boundary conditions is $b_z(x) = z -z \left({\frac{z-1}{z}}\right)^x$. Via Eq.\eqref{eq:resolvent} the resolvent becomes
\begin{equation}
	G(z)=\int_0^1\!dx\, z^{x-1}(z-1)^{-x}
\end{equation}
and has a branch cut at $z\in[0,1]$. Cauchy's identity yields the spectral density as $G(\lambda-i\epsilon)-  G(\lambda+i\epsilon)=2i\pi\sigma_{[0,1]}(\lambda)$ from which we find $d\sigma_{[0,1]}(\lambda)=  \frac{d\lambda}{\pi}\int_0^1dx\,\sin(\pi x)\, \lambda^{x-1}(1-\lambda)^{-x}$. Integrating over $x$ leads to
\begin{equation}\label{eq:spectrum_01}
	d\sigma_{[0,1]}(\lambda) = \frac{d\lambda}{\lambda(1-\lambda)}\,
	\frac{1}{\pi^2 + \log^2(\frac{1-\lambda}{\lambda})}1_{\lambda\in[0,1]}.
\end{equation}
By a change of variable, this is actually a Cauchy-Lorenz distribution. Defining $\nu:= \log(\frac{\lambda}{1-\lambda})\in(-\infty,+\infty)$, or $\lambda=\frac{e^{\nu}}{1+e^{\nu}}$, we have
\[
d\sigma_{[0,1]}(\lambda) =  \frac{ d\nu}{\pi^2+\nu^2}.
\]

\begin{figure}[h]
	\centering\includegraphics[width=0.7 \linewidth]{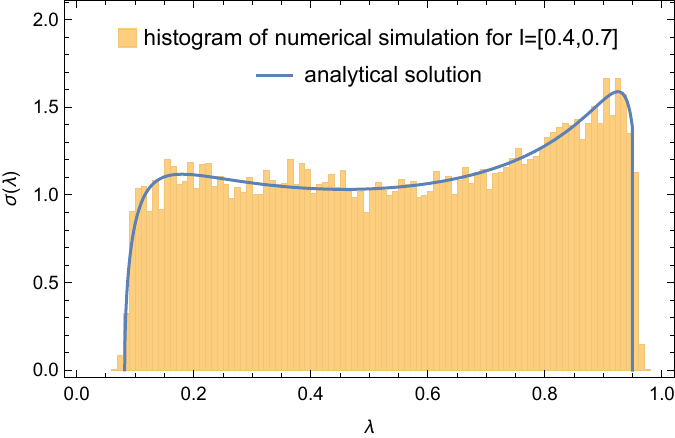}
	\caption{\label{fig:numerical_sim}Comparison between the analytical solution and a numerical simulation for the spectrum $\sigma_I$ of $M_I$ for $I=[0.4,0.7]$. The numerics comes from a simulation of the QSSEP time evolution of $M$ in Eq.\eqref{eq:M_evolution} on N = 100 sites with discretization time step $dt = 0.1$. Instead of averaging over many noisy realizations, we exploit the ergodicity of QSSEP and perform a time average over a single realization between t = 0.25 and t = 0.4. The QSSEP dynamics reaches its steady state at approximately t = 0.2.}
\end{figure}

\paragraph{Spectrum of a subblock of $M$.}
The spectrum of an arbitrary subblock $M_I\subset M$ can become quite complicated as the following proposition shows. Since we are dealing with a structured matrix, the spectrum will depends on the position of the subblock and not on its size only. Figure~\ref{fig:numerical_sim} shows that the analytical result from the following proposition indeed agrees with a numerical simulation of the spectrum of $M_I$.

\medskip
\begin{proposition}
	The spectrum of subblock $M_I$ of $M$ restricted to the interval $I=[c,d]$ of size $\ell=d-c$ is
	\begin{equation}
		d\sigma_{[c,d]}(\lambda)=   \frac{d\lambda}{\pi\lambda(1-\lambda)}\,\frac{\theta}{\theta^2 + (\log r)^2}\, 1_{\lambda\in[z_*^-,z_*^+]}
	\end{equation}
	where $r$ and $\theta$ are functions of $\lambda$ determined by the complex transcendental equation
	\begin{equation}
		\frac{\lambda}{1-\lambda} re^{i\theta} = \frac{c(\log r+ i\theta)  -\ell }{(1-d)(\log r  + i\theta) +\ell}.
	\end{equation}
	The support of the spectrum is $[z_*^-,z_*^+]$ with
	\begin{equation} \label{eq:z*-pm}
		z_*^\pm = \frac{c(1-c)+d(1-d) \pm \sqrt{\Delta}}{c(1-c)+d(1-d) \pm \sqrt{\Delta} + 2(1-d)^2\,e^{-\delta_\pm}}.
	\end{equation}
	with $\Delta=\ell(1-\ell)[\ell(1-\ell)+4c(1-d)]$ and $\ell+c\delta^\pm= \frac{1}{2(1-d)}[\ell(1-\ell)\pm\sqrt{\Delta}]$.
	Note that  $z_*^-<c$ and $d<z_*^+$ so that the support of the spectrum is larger than the interval $[c,d]$.
\end{proposition}

\begin{proof}
An ansatz for the differential equation Eq.\eqref{eq:b_z_equation}, satisfying $b_z(0)=0$ and $b_z(1)=1$, is
\begin{equation} 
	 b_z(x) = 
	 \begin{cases}
	 	\alpha x, & \text{if $0<x<c$}, \\
	 	z + \gamma\, Q(z)^y, & \text{if $c<x<d$}, \\ 
	 	1+\beta(x-1), & \text{if $d<x<1$},
	 \end{cases}
\end{equation}
where instead of $x\in[c,d]$ we use $y\in[0,1]$ to parametrize the interval, $x=c+y\ell$. The complex function $Q(z)$ parametrises the exponential growth of $z-b_z(x)$ in the interval $[c,d]$. The coefficients $\alpha,\, \beta,\, \gamma$ are determined, as a function of $Q$, by the continuity of $b_z$ and $b_z'$ at the boundaries of the interval $y=0$ and $y=1$. One finds (we don't need the explicit expression for $\alpha$ and $\beta$)
\begin{equation} \label{eq:eqforB}
b_z(x) = z + \frac{c(1-z) - (1-d) zQ(z)}{(1 -\ell)Q(z)}\, Q(z)^y, 
\end{equation}
for $x=c+y\ell\in[c,d]$. The continuity of $b_z$ and $b_z'$ at the boundaries of the interval yields four equations. Three are used to determined $\alpha,\, \beta, \gamma$. The fourth one yields a constraint on $Q$,
\begin{equation} \label{eq:eqfor Q}
(1-z+ zQ)(\ell-c\log Q) = z(\ell -1) Q\log Q.
\end{equation}
This specifies the analytical structure of $Q(z)$, as complex function of $z$, from which we deduce the spectral density.

We first determine the support of the spectrum -- by looking at the position of the cut of the function $Q(z)$ -- and then the spectral measure on that support -- by computing the jump of $Q(z)$ on its cut.

The cut of $Q$ is on the real axis. To find it, we write Eq.\eqref{eq:eqfor Q} as, 
\[  
f(\hat Q)=g(\hat Q),\quad \mathrm{with}\ f(q) := 1 +\eta q,\ g(q):= \frac{(1-\ell)\log q}{\ell + c\log q},
\]
where $\hat Q:=Q^{-1}$ and $\eta := \frac{1-z}{z}$.
We have $g'(q)= \frac{\ell(1-\ell)}{q(\ell+c\log q)^2}>0$, and $g(q)$ diverges for $q=e^{-\ell/c}<1$ and $g(0)=g(\infty)=\frac{1-\ell}{c}>1$.
On $\mathbb{R}_+$, the function $g(q)$ is critical to the straight-line $f(q)$ for two critical values $\eta_*^\pm$ with $\eta_*^+<\eta_*^-$. 
This corresponds to two critical values for $z$, i.e. $z_*^\pm= 1/(1+\eta_*^\pm)$ with $0< z_*^-< z_*^+<1$.
The cut, and hence the support of the eigenvalues, is thus on $[z_*^{-},z_*^{+}]\subset[0,1]$. 

The two critical values for $\eta$ are solutions of (with $q=e^\delta$)
\[ 
1+\eta e^{\delta} =\frac{ (1-\ell)\delta}{\ell+c\delta},\ \eta e^{\delta} = \frac{ \ell(1-\ell)}{(\ell+c\delta)^2} .
\]
This leads to a second order equation for $\delta$, 
\[
c(1-d) \delta^2 + \ell(1-d-c)\delta - \ell=0
\]
The discriminant is $\Delta= \ell(1-\ell)[\ell(1-\ell)+4c(1-d)]$ and the two solutions $\delta^\pm := \frac{1}{2c(1-d)}[\ell(d+c-1)\pm\sqrt{\Delta}]$. We let $\kappa_\pm:= \ell(1-\ell)\pm\sqrt{\Delta}$, which is symmetric by $[c,d]\to [1-d,1-c]$. Then, $\eta_\pm=4(1-d)^2\ell(1-\ell) e^{-\delta_\pm}/\kappa_\pm^2$, so that
\begin{equation}
	z_*^\pm = \frac{c(1-c)+d(1-d) \pm \sqrt{\Delta}}{c(1-c)+d(1-d) \pm \sqrt{\Delta} + 2(1-d)^2e^{-\delta_\pm}}.
\end{equation}
as proposed in Eq.\eqref{eq:z*-pm}.

Let us now compute the spectral density. The later is obtained by integrating the branch cut discontinuity of the resolvent \eqref{eq:resolvent} (with the pole at the origin discarded), so that
\[
d\sigma_{[c,d]}(\lambda)= \Im{m} \, \frac{d\lambda}{\pi} \int_c^d \frac{dx}{\ell} \frac{1}{z- b_z(x)}
\]
with $z-b_z(x)$ given by Eq.\eqref{eq:eqforB} for $x=c+y\ell$ in $[c,d]$.
Recall Eq.\eqref{eq:eqfor Q} for $Q$, which can alternatively be written as $(c(1-z)-(1-d)zQ)\log Q = \ell (1-z+zQ)$ so that 
\[
z-  b_z(x) = - (\frac{\ell}{1-\ell})( \frac{1-z+zQ}{Q})\, \frac{Q^y}{\log Q},
\]
Using $Q^{-y}\log Q= -\partial_y Q^{-y}$, we can explicitly integrate the discontinuity to get
\begin{equation} \label{eq:Sigma1rst}
d\sigma_{[c,d]}(\lambda)=  \frac{d\lambda}{\pi}  (\frac{1-\ell}{\ell})\, \Im{m} \left[\frac{1-Q}{1-\lambda+\lambda Q} \right]\,1_{\lambda\in[z_*^-,z_*^+]} 
\end{equation}
Eq.\eqref{eq:eqfor Q} also allows to express $Q$ as a function of $\log Q$ as
\begin{equation} \label{eq:2ndeqforQ}
	\frac{z}{1-z}Q = \frac{ c\log Q - \ell}{(1-d)\log Q + \ell}.
\end{equation}
This can then be used to simplify further the expression \eqref{eq:Sigma1rst} for the spectral density as
\begin{equation} \label{eq:Sigma_final}
d\sigma_{[c,d]}(\lambda)=  \frac{d\lambda}{\pi \lambda(1-\lambda) } \Im{m} \left[\frac{1}{\log Q}\right]\, 1_{\lambda\in[z_*^-,z_*^+]}
\end{equation}
Parametrising $Q$ as $Q=re^{i\theta}$ in Eqs.(\ref{eq:2ndeqforQ},\ref{eq:Sigma_final}) yields the claim.
\end{proof}

To end this section, let us note that the explicit expression of the spectral density satisfies the expected symmetry $d\sigma_{[c,d]}(\lambda)= d\sigma_{[1-d,1-c]}(1-\lambda)$. In particular, one can verify the symmetry of the support, that is $z_*^\pm(1-d, 1-c)= 1- z_*^\mp(c,d)$. Furthermore, for $d=1^-$, $\ell=1-c+0^+$, we have $\delta_+=+\infty$ and $\delta_-=-1/c$ so that $z_*^+(c,1)=1$  and $z_*^-(c,1)= \frac{c}{c+(1-c)e^{-1/c}}$. By symmetry, for $c=0^+$, $\ell=1-d-0^+$, we have $\delta_-=-\infty$ and $\delta_+=1/(1-d)$ so that $z_*^-(0,d)=0$ and $z_*^+(0,d)=\frac{d}{d + (1-d)e^{-1/(1-d)}}$. For $c=0^+$, $d=1^-$, we get $z_*^-=0$, $z_*^+=1$ and $\theta=\pi$, $r=\frac{1-\lambda}{\lambda}$, and we recover Eq.\eqref{eq:spectrum_01}.

%In both cases, we recover the known results \cite{Bernard2023Exact}. For $c=0^+$, $d=1^-$, we recover \eqref{eq:spectrum_01}.

\bmhead{Acknowledgments}
We thank Philippe Biane for discussions and Roland Speicher for suggesting to look at non-linear operations on these matrix ensembles.
This work was in part supported by the CNRS, the ENS and the ANR project ESQuisses under contract number ANR-20-CE47-0014-01.

\section*{Declarations}

\bmhead{Conflict of interest}
	On behalf of all authors, the corresponding author states that there is no conflict of
	interest.
\bmhead{Availability of data}
	Data sharing not applicable to this article as no external datasets were analysed during the current study

\begin{appendices}

\section{Free probability glossary}\label{app:glossary}
We tried to use notation as close as possible to those used by R. Speicher in \cite{Speicher2019Lecture}.
Let $\sigma$ be a (classical) measure for a random variable $X$.
\begin{itemize}
	\item 
The resolvent $G_\sigma(z):=\mathbb{E}_\sigma[\frac{1}{z-X}]$ is a generating function of the moments $m_n:=\mathbb{E}_\sigma[X^n]$:
\begin{align*}
	G_\sigma(z) &:= \sum_{n\geq 0} m_n\, z^{-n-1} = z^{-1} + m_1 z^{-2} + m_2 z^{-3} +\cdots\\
	\hat M_\sigma(z)&:= z^{-1}G_\sigma(z^{-1}) = \sum_{n\geq 0} m_n\, z^{n} =  1 + m_1 z + m_2 z^{2} +\cdots
\end{align*}
(We put a hat on $\hat M$ to distinguished from the matrix $M$, otherwise it is the same definition as in Speicher et al).

\item 
The $R$-transform is a generating function of the free cumulants $\kappa_p:=\kappa_p(\sigma)$:
\begin{align*}
	R_\sigma(z) &:= \sum_{p\geq 1} \kappa_p\, z^{p-1} = \kappa_1  + \kappa_2 z + \kappa_3 z^2 +\cdots \\
	K_\sigma(z) &:=z^{-1}+R_\sigma(z) = \sum_{p\geq 0} \kappa_p\, z^{p-1} = z^{-1}+ \kappa_1  + \kappa_2 z + \kappa_3 z^2 +\cdots \\
	C_\sigma(z) &:=zR_\sigma(z)= \sum_{p\geq 1} \kappa_p\, z^{p} = \kappa_1 z  + \kappa_2 z^2 + \kappa_3 z^3 +\cdots
\end{align*} 
(There is a shift of $1$ in this definition of $C_\sigma$ compared to that of Speicher et al).
We of course have $zK_\sigma(z)=1+C_\sigma(z)$.

\item
The function $G_\sigma$ and $K_\sigma$ are inverse each other, thus
\begin{equation*}
	K_\sigma(G_\sigma(z))=z,\quad G_\sigma(K_\sigma(z))=z
\end{equation*}
The previous relation then reads
\begin{equation*}
	zG_\sigma(z) = 1 + C_\sigma(G_\sigma(z))\ ,\ 
	\hat M_\sigma(z) = 1 + C_\sigma(z\hat M_\sigma(z)) .
\end{equation*}
\item 
The $S$-transform can be defined by
\begin{equation*}
	C_\sigma(zS_\sigma(z))=z,\quad C_\sigma(z)\,S_\sigma(C_\sigma(z))=z
\end{equation*}
The function $S_\sigma$ exists, as a formal power in $z$, whenever $\kappa_1\not= 0$ : $S_\sigma(z) = \frac{1}{\kappa_1} - \frac{\kappa_2}{\kappa_1^3}z + \cdots$.
Using $G_\sigma(K_\sigma(z))=z$, this relation can alternatively be written as
\[
G_\sigma\left(\frac{1+z}{zS_\sigma(z)}\right)= zS_\sigma(z),\quad S_\sigma(zG_\sigma(z)-1)=\frac{G_\sigma(z)}{zG_\sigma(z)-1} .
\]
Setting $w=zG_\sigma(z)-1$, the above formula can be written as $S_\sigma(w)=\frac{w+1}{zw}$ with $z(w)$ determined by solving $zG_\sigma(z)=w+1$.

\item 
For two measures $\sigma$ and $\nu$, the additive free convolution is defined 
\[ R_{\sigma \boxplus\nu}(z) = R_{\sigma}(z)  + R_{\nu}(z) ,\]
that is, we add the free cumulants. Thus if $a$ and $b$ are (relatively) free then $R_{a+b}(z)=R_a(z)+R_b(z)$.

\item 
For two measures $\sigma$ and $\nu$, the free multiplicative convolution $\sigma \boxtimes\nu$ is defined via their $S$-transform
\[ S_{\sigma \boxtimes\nu}(z) = S_{\sigma}(z)  S_{\nu}(z) ,\]
that is, we multiply the $S$-transforms. Thus, if $a$ and $b$ are (relatively) free, then $S_{ab}(z)=S_a(z)S_b(z)$ (instead of $ab$ we could have considered $a^{\frac{1}{2}}ba^{\frac{1}{2}}$).
\end{itemize}

%\medskip

%\item 
%$R$ \& $S$ transform of bi-modal distribution.\\
%We take $\sigma=\ell \delta_1 + (1-\ell)\delta_0$ as measure, so that the resolvent is $G(z)= \frac{1-\ell}{z}+\frac{\ell}{z-1}$. This corresponds to the spectral measure of a projector, with $\ell$ the rank of the projector.
%
%The $R$-transform is such that $G(K(z))=z=K(G(z))$ so that is satisfies $zK(K-1) = \ell K + (1-\ell)(K-1)$. Hence
%\[
%K(z) = \frac{1}{2z}(z+1 + \sqrt{(z-1)^2 + 4z\ell})  .
%\]
%We checked that $K(z)=\frac{1}{z}+\cdots$ as $z\to 0$. 
%One defines $C(z)$ via $zK(z)=1+C(z)$, so that $C(z) = \frac{1}{2}(z-1 + \sqrt{(z-1)^2 + 4z\ell})$, or 
%\[
%R(z) = \frac{1}{2z}(z-1 + \sqrt{(z-1)^2 + 4z\ell})  .
%\]
%
%The $S$-transform is defined by $C(zS(z))=z$ so that it satisfies $z+1-zS=\ell S$ and hence
%\[
%S(z)=\frac{z+1}{z+\ell} .
%\]
%We check that $S(z)=\frac{1}{\ell}+\cdots$ as $z\to 0$ (and $\ell$ is the first moment of the bi-modal distribution).

\section{Local free cumulants for Haar-randomly rotated matrices}\label{app:free_cum_haar}
Let $M=UDU^\dag$, with $U$ Haar distributed over the unitary group and $D$ a diagonal matrix with spectral density $\sigma$ in the large $N$ limit.
From the HCIZ integral, it is known that the generating function $\mathbb{E}[e^{N\tr(QM)}]$ can be written in terms of the free cumulants $\kappa_n(\sigma)$ of the density $\sigma$ as \cite{Zinn-Justin_2003}
\[
\mathbb{E}[e^{z N \tr(\hat QM)}] \asymp_{N\to\infty} \exp\Big( N \sum_{k\geq1} \frac{z^k}{k}\tr(\hat Q^k)\,\kappa_k(\sigma) \Big)
\]
for any finite rank matrix $\hat Q$. 

Let us prove (see also \cite{Maillard_2019}) that this implies that the local free cumulants are $g_n=\kappa_n(\sigma)$, that is
\begin{align} \label{eq:Haar-loop}
	\mathbb{E}[M_{12}M_{23}\cdots M_{n1}]=N^{1-n}\, \kappa_n(\sigma)\, (1+O(N^{-1}))
\end{align}
Note that due to $U(N)$ invariance (which in particular includes permutations), all sets of distinct indices $i_1,i_2,\cdots,i_n$ are equivalent.

Choose $\hat Q=P_n$ the cyclic permutation $(12\cdots n)$, so that $\tr(P_nM)=M_{12}+M_{23}+\cdots +M_{n1}$. It is easy to see (using $U(1)^N\subset U(N)$ invariance), that the first non-vanishing term in $\mathbb{E}[e^{zN \tr(P_nM)}]$ is of order $z^n$ and given by $z^n N^n \mathbb{E}[(\tr(P_nM))^n]$. Furthermore, (this can be proved say by recurrence)
\[
\mathbb{E}[(\tr(P_nM))^n] = \mathbb{E}[(M_{12}+M_{23}+\cdots +M_{n1})^n]{=} n!\, \mathbb{E}[M_{12}M_{23}\cdots M_{n1}]
\]
Thus 
\[
\mathbb{E}[e^{z N \tr(P_nM)}] = z^nN^n\, \mathbb{E}[M_{12}M_{23}\cdots M_{n1}] + O(z^{n+1})
\]
Since $\tr(P_n^k)=0$ for $k<n$ and $\tr(P_n^n)=n$, we have 
\[
e^{N \sum_{k\geq0} \frac{z^k}{k}\tr(P_n^k)\,\kappa_k(\sigma)}= N z^n \kappa_n(\sigma)+ O(z^{n+1})
\]
Comparing the two last equations proves Eq.\eqref{eq:Haar-loop}.

\end{appendices}


\bibliography{references}

\end{document}